\documentclass[11pt]{article}
\usepackage{amsfonts}
\usepackage{amsmath}
\usepackage{graphicx}
\usepackage{url}

\newcommand{\eq}[1]{(\ref{eq:#1})}

\newtheorem{theorem}{Theorem}[section]
\newtheorem{definition}[theorem]{Definition}
\newtheorem{lemma}[theorem]{Lemma}
\newtheorem{proposition}[theorem]{Proposition}
\newtheorem{corollary}[theorem]{Corollary}

\numberwithin{equation}{section}

\newcommand{\qed}{\rule{7pt}{7pt}}
\newenvironment{proof}{\noindent{\bf Proof}\hspace*{1em}}{\hfill\qed\vspace{0.125in}}

\newcommand{\R}{\mathbb{R}}
\newcommand{\bfA}{\mathbf{A}}
\newcommand{\cA}{\mathcal{A}}
\newcommand{\id}{I}

\newcommand{\trace}{\operatorname{Tr}}
\newcommand{\diag}{\operatorname{diag}}
\newcommand{\sign}{\operatorname{sign}}
\newcommand{\rank}[1]{\operatorname{rank}(#1)}
\newcommand{\card}[1]{\operatorname{card}(#1)}
\newcommand{\vvec}[1]{\operatorname{vec}(#1)}

\begin{document}

\title{Guaranteed Minimum-Rank Solutions of Linear Matrix Equations via Nuclear Norm Minimization}
\author{Benjamin Recht \thanks{Center for the Mathematics of Information, California Institute of Technology} \and
Maryam Fazel \thanks{Control and Dynamical Systems, California Institute of Technology}
\and
Pablo A. Parrilo \thanks{Laboratory for Information and Decision Systems, Massachusetts Institute of Technology}}
\date{\today}
\maketitle

\begin{abstract}
The affine rank minimization problem consists of finding a matrix of
minimum rank that satisfies a given system of linear equality
constraints. Such problems have appeared in the literature of a
diverse set of fields including system identification and control,
Euclidean embedding, and collaborative filtering. Although specific
instances can often be solved with specialized algorithms, the
general affine rank minimization problem is NP-hard, because it
contains vector cardinality minimization as a special case.

In this paper, we show that if a certain restricted isometry
property holds for the linear transformation defining the
constraints, the minimum rank solution can be recovered by solving a
convex optimization problem, namely the minimization of the nuclear
norm over the given affine space. We present several random
ensembles of equations where the restricted isometry property holds
with overwhelming probability, provided the codimension of the
subspace is $\Omega(r (m+n) \log mn)$, where $m,n$ are the
dimensions of the matrix, and $r$ is its rank.

The techniques used in our analysis have strong parallels in the
compressed sensing framework. We discuss how affine rank
minimization generalizes this pre-existing concept and outline a
dictionary relating concepts from cardinality minimization to those
of rank minimization. We also discuss several algorithmic approaches
to solving the norm minimization relaxations, and illustrate our
results with numerical examples.
\end{abstract}

\noindent{\footnotesize {\bf Keywords.} rank, convex optimization,
matrix norms, random matrices, compressed sensing, semidefinite
programming.}

\section{Introduction}

Notions such as order, complexity, or dimensionality can often be
expressed by means of the rank of an appropriate matrix. For
example, a low-rank matrix could correspond to a low-degree
statistical model for a random process (e.g., factor analysis), a
low-order realization of a linear system~\cite{Fazel01}, a low-order
controller for a plant~\cite{ElGhaoui93}, or a low-dimensional
embedding of data in Euclidean space~\cite{Linial95}. If the set of
feasible models or designs is affine in the matrix variable,
choosing the simplest model can be cast as an \emph{affine rank
minimization problem},
\begin{equation}
\begin{array}{ll}
\mbox{minimize} & \rank{X}\\
\mbox{subject to} & \cA(X)=b,
\end{array}
\label{eq:min-rank-prob}
\end{equation}
where $X\in \R^{m \times n}$ is the decision variable, and the
linear map $\cA:\R^{m \times n} \rightarrow \R^p$ and vector
$b \in \R^p$  are given.  In certain instances with very special
structure, the rank minimization problem can be solved by using the
singular value decomposition, or can be exactly reduced to the
solution of linear systems~\cite{Mesbahi97,PKrank}. In general,
however, problem~\eq{min-rank-prob} is a challenging nonconvex
optimization problem for which all known finite time algorithms have
at least doubly exponential running times in both theory and
practice. For the general case, a variety of heuristic algorithms
based on local optimization, including alternating
projections~\cite{Grigoriadis00} and alternating
LMIs~\cite{Skelton98}, have been proposed.

A recent heuristic introduced in~\cite{FazelThesis} minimizes the
\emph{nuclear norm}, or the sum of the singular values of the
matrix, over the affine subset.  The nuclear norm is a convex
function, can be optimized efficiently, and is the best convex
approximation of the rank function over the unit ball of matrices
with norm less than one.  When the matrix variable is symmetric and
positive semidefinite, this heuristic is equivalent to the trace
heuristic often used by the control community (see,
e.g.,~\cite{Beck98,Mesbahi97}). The nuclear norm heuristic has been
observed to produce very low-rank solutions in practice, but a
theoretical characterization of when it produces the minimum rank
solution has not been previously available.  This paper provides the
first such mathematical characterization.

Our work is built upon a large body of literature on a related
optimization problem.  When the matrix variable is constrained to be
diagonal, the affine rank minimization problem reduces to finding
the \emph{sparsest vector} in an affine subspace. This problem is
commonly referred to as \emph{cardinality minimization}, since we
seek the vector whose support has the smallest cardinality, and is
known to be NP-hard~\cite{Natarajan95}.  For diagonal matrices, the
sum of the singular values is equal to the sum of the absolute
values (i.e., the $\ell_1$ norm) of the diagonal elements.
Minimization of the $\ell_1$ norm is a well-known heuristic for the
cardinality minimization problem, and stunning results pioneered by
Cand\`es and Tao~\cite{CandesCS} and Donoho~\cite{DonohoCS} have
characterized a vast set of instances for which the $\ell_1$
heuristic can be \emph{a priori} guaranteed to yield the optimal
solution. These techniques provide the foundations of the recently
developed \emph{compressed sensing} or \emph{compressive sampling}
frameworks for measurement, coding, and signal estimation. As has
been shown by a number of research groups (e.g.,
\cite{Wakin07,Candes05b,CRT06,Candes05}), the $\ell_1$ heuristic for
cardinality minimization provably recovers the sparsest solution
whenever the sensing matrix has certain ``basis incoherence''
properties, and in particular, when it is randomly chosen according
to certain specific ensembles.

The fact that the $\ell_1$ heuristic is a special case of the
nuclear norm heuristic suggests that these results from the
compressed sensing literature might be extended to provide
guarantees about the nuclear norm heuristic for the more general
rank minimization problem. In this paper, we show that this is
indeed the case, and the parallels are surprisingly strong.
Following the program laid out in the work of Cand\`es and Tao, our
main contribution is the development of a restricted isometry
property (RIP), under which the nuclear norm heuristic can be
\emph{guaranteed} to produce the minimum-rank solution. Furthermore,
as in the case for the $\ell_1$ heuristic, we provide several
specific examples of matrix ensembles for which RIP holds with
overwhelming probability. Our results considerably extend the
compressed sensing machinery in a so far undeveloped direction, by
allowing a much more general notion of parsimonious models that rely
on low-rank assumptions instead of cardinality restrictions.

To make the parallels as clear as possible, we begin by establishing
a dictionary between the matrix rank and nuclear norm minimization
problems and the vector sparsity and $\ell_1$ norm problems in
Section~\ref{sec:dictionary}.  In the process of this discussion, we
present a review of many useful properties of the matrices and
matrix norms necessary for the main results. We then generalize the
notion of Restricted Isometry to matrices in Section~\ref{sec:rip}
and show that when linear mappings are Restricted Isometries,
recovering low-rank solutions of underdetermined systems can be
achieved by nuclear norm minimization. In Section
\ref{sec:gaussian}, we present several families of random linear
maps that are restricted isometries with overwhelming probability
when the dimensions are sufficiently large. In
Section~\ref{sec:algorithms}, we briefly discuss three different
algorithms designed for solving the nuclear norm minimization
problem and their relative strengths and weaknesses: interior point
methods, gradient projection methods, and a low-rank factorization
technique. In Section~\ref{sec:experiments}, we demonstrate that in
practice nuclear-norm minimization recovers the lowest rank
solutions of affine sets with even fewer constraints than those
guaranteed by our mathematical analysis.  Finally, in
Section~\ref{sec:conclusions}, we list a number of possible
directions for future research.

\subsection{When are random constraints interesting for rank minimization?}
\label{sec:applications}

As in the case of compressed sensing, the conditions we derive to
guarantee properties about the nuclear norm heuristic are
deterministic, but they are at least as difficult to check as
solving the rank minimization problem itself. We are only able to
guarantee that the nuclear norm heuristic recovers the minimum rank
solution of $\cA(X)=b$ when $\cA$ is sampled from
specific ensembles of random maps.  The constraints appearing in
many of the applications mentioned above, such as low-order control
system design, are typically not random at all and have structured
demands according to the specifics of the design problem. It thus
behooves us to present several examples where random constraints
manifest themselves in practical scenarios for which no practical
solution procedure is known.

\paragraph{Minimum order linear system realization}
Rank minimization forms the basis of many model reduction and
low-order system identification problems for linear time-invariant
(LTI) systems. The following example illustrates how random
constraints might arise in this context. Consider the problem of
finding the minimum order discrete-time LTI system that is
consistent with a set of time-domain observations. In particular,
suppose our observations are the system output sampled at a fixed
time $N$, after a random Gaussian input signal is applied from $t=0$
to $t=N$. Suppose we make such measurements for $p$ different input
signals, that is, we observe $y_i(N)=\sum_{t=0}^{N} a_i(N-t)h(t)$
for $i=1,\ldots,p$, where $a_i$,  the $i$th input signal, is a
zero-mean Gaussian random variable with the same variance for
$t=0,\ldots N$, and $h(t)$ denotes the impulse response. We can
write this compactly as $y=Ah$, where $h=[h(0),\ldots,h(N)]'$, and
$A_{ij}=a_i(N-j)$.

From linear system theory, the order of the minimal realization for
such a system is given by the rank of the following Hankel matrix
(see, e.g.,~\cite{Fazel03,Sontag98})
\[
\operatorname{hank}(h) :=
  \left[\begin{array}{cccc}h(0)   & h(1)    & \cdots & h(N)\\
                h(1)   & h(2)  & \cdots & h(N+1)\\
               \vdots & \vdots &        &\vdots\\
                h(N)   & h(N+1) & \cdots & h(2N)\\ \end{array}\right].
\]
Therefore the problem can be expressed as
\[
\begin{array}{ll}
\mbox{minimize} & \rank{\operatorname{hank}(h)}\\
\mbox{subject to} & A h=y
\end{array}
\]
where the optimization variables are $h(0),\ldots,h(2N)$, and the
matrix $A$ consists of i.i.d. zero-mean Gaussian entries.

\paragraph{Low-Rank Matrix Completion}

In the matrix completion problem where we are given  random subset
of entries of a matrix, we would like to fill in the missing entries
such that the resulting matrix has the lowest possible rank.  This
problem arises in machine learning scenarios where we are given
partially observed examples of a process with a low-rank covariance
matrix and would like to estimate the missing data.  A typical
situation where the hidden matrix is low-rank is when the columns
are i.i.d. samples of a random process with low-rank covariance.
Such models are ubiquitous in Factor Analysis, Collaborative
Filtering, and Latent Semantic
Indexing~\cite{Rennie05,SrebroThesis}.  In many of these settings,
some prior probability distribution (such as a Bernoulli model or
uniform distribution on subsets) is assumed to generate the set of
available entries.

Suppose we are presented with a set of triples $(I(i),J(i),S(i))$
for $i=1,\ldots,p$ and wish to find the matrix with $S(i)$ in the
entry corresponding to row $I(i)$ and column $J(i)$ for all $i$. The
matrix completion problem seeks to solve
\begin{equation*}
    \begin{split}
        \min_{Y} &\quad \rank{Y}\\
        \mbox{s.t.} &\quad Y_{I(i),J(i)} = S(i), \qquad i=1,\ldots, K
    \end{split}
\end{equation*}
which is a special case of the affine rank minimization problem.

\paragraph{Low-dimensional Euclidean embedding problems}

A problem that arises in a variety of fields is the determination of
configurations of points in low-dimensional Euclidean spaces,
subject to some given distance information. In Multi-Dimensional
Scaling (MDS), such problems occur in extracting the underlying
geometric structure of distance data. In psychometrics, the
information about inter-point distances is usually gathered through
a set of experiments where subjects are asked to make quantitative
(in metric MDS) or qualitative (in non-metric MDS) comparisons of
objects. In computational chemistry, they come up in inferring the
three-dimensional structure of a molecule (molecular conformation)
from information about interatomic distances~\cite{Trosset00}.

A symmetric matrix $D\in \mathcal{S}^{n}$ is called a \emph{Euclidean
distance matrix} (EDM) if there exist points $x_1, \ldots, x_n$ in
$\R^d$ such that $D_{ij}=\|x_i-x_j\|^2$. Let
$V:=I_n-\frac{1}{n}\mathbf{1}\mathbf{1}^T$ be the projection matrix
onto the hyperplane $\{v\in\R^n\,\,:\,\,\mathbf{1}^Tv=0\}$.  A
classical result by Schoenberg states that $D$ is a Euclidean distance
matrix of $n$ points in $\R^d$ if and only if $D_{ii}=0$, the matrix
$VDV$ is negative semidefinite, and $\rank{VDV}$ is less than or equal
to $d$~\cite{Schoenberg35}. If the matrix $D$ is known exactly, the
corresponding configuration of points (up to a unitary transform) is
obtained by simply taking a matrix square root of $-\frac{1}{2}VDV$.
However, in many cases, only a random sampling collection of the
distances are available. The problem of finding a valid EDM consistent
with the known inter-point distances and with the smallest embedding
dimension can be expressed as the rank optimization problem
\begin{equation*}
\label{eq:edm-rank}
\begin{array}{lll}
\mbox{minimize} & \rank{VDV}\\
\mbox{subject to} &  VDV \preceq 0\\
& \cA(D)=b,
\end{array}
\end{equation*}
where $\cA:\mathcal{S}^{n} \rightarrow \R^{p}$ is a random
sampling operator as discussed in the matrix completion problem.

This problem involves a Linear Matrix Inequality (LMI) and appears
to be more general than the equality constrained rank minimization
problem. However, general LMIs can equivalently be expressed as rank
constraints on an appropriately defined block matrix. The rank of a
block symmetric matrix is equal to the rank of a diagonal block plus
the rank of its Schur complement (see, e.g.,~\cite[\S
2.2]{HornJohnsonBook2}). Given a function $f$ that maps matrices into
$q\times q$ symmetric matrices, that $f(X)$ is positive semidefinite
can be equivalently expressed through a rank constraint as
\[
f(X) \succeq 0 \qquad \Leftrightarrow \qquad
\operatorname{rank} \left(\left[\begin{array}{cc} \id_q & B \\
             B' & f(X) \end{array}\right]\right) \leq q, \quad \mbox{ for some }
B \in \R^{q \times q}.
\]
That is, if there exists a matrix $B$ satisfying the inequality
above, then $f(X)=B'B\succeq 0$. Using this equivalent
representation allows us to rewrite \eq{edm-rank} as an affine rank
minimization problem.

\paragraph{Image Compression}

A simple and well-known method to compress two-dimensional images
can be obtained by using the singular value decomposition
(e.g.,~\cite{imageSVD}). The basic idea is to associate to the given
grayscale image a rectangular matrix $M$, with the entries $M_{ij}$
corresponding to the gray level of the $(i,j)$ pixel. The best
rank-$k$ approximation of $M$ is given by
\[
X^* := \arg \min_{\rank{X} \leq k} ||M-X||,
\]
where $||\cdot||$ is any unitarily invariant norm.  By the classical
Eckart-Young-Mirsky theorem (\cite{EckartYoung, Mirsky}), the
optimal approximant is given by a truncated singular value
decomposition of $M$, i.e., if $M = U \Sigma V^T$, then $X^* = U
\Sigma_k V^T$, where the first $k$ diagonal entries of $\Sigma_k$
are the largest $k$ singular values, and the rest of the entries are
zero.  If for a given rank $k$, the approximation error $||M-X^*||$
is small enough, then the amount of data needed to encode the
information about the image is $k (m+n-k)$ real numbers, which can
be much smaller than the $mn$ required to transmit the values of all
the entries.

Consider a given image, whose associated matrix $M$ has low-rank, or
can be well-approximated by a low-rank matrix. As proposed by Wakin
\emph{et al.}\ \cite{singlepixel}, a single-pixel camera would
ideally produce measurements that are random linear combinations of
all the pixels of the given image. Under this situation, the image
reconstruction problem boils down exactly to affine rank
minimization, where the constraints are given by the random linear
functionals.

It should be remarked that the simple SVD image compression scheme
described has certain deficiencies that more sophisticated
techniques do not share (in particular, the lack of invariance of
the description length under rotations). Nevertheless, due to its
simplicity and relatively good practical performance, this method is
particularly popular in introductory treatments and numerical linear
algebra textbooks.

\section{From Compressed Sensing to Rank Minimization}\label{sec:dictionary}

As discussed above, when the matrix variable is constrained to be
diagonal, the affine rank minimization problem~\eq{min-rank-prob}
reduces to the cardinality minimization problem of finding the
element in the affine space with the fewest number of nonzero
components. In this section we will establish a dictionary between
the concepts of rank and cardinality minimization.  The main
elements of this correspondence are outlined in
Table~\ref{table:dictionary}. With these elements in place, the
existing proofs of sparsity recovery provide a template for the more
general case of low-rank recovery.

In establishing our dictionary, we will provide a review of useful
facts regarding matrix norms and their characterization as convex
optimization problems.  We will show how computing both the operator norm and
the nuclear norm of a matrix can be cast as semidefinite programming
problems. We also establish the suitable optimality conditions for
the minimization of the nuclear norm under affine equality
constraints, the main convex optimization problem studied in this
article. Our discussion of matrix norms will mostly follow the
discussion in~\cite{FazelThesis,Vandenberghe96} where extensive
lists of references are provided.

\begin{table}
\begin{center}
\begin{tabular}{|c||c|c|}
\hline
{\bf  parsimony concept} & cardinality & rank \\
\hline {\bf Hilbert Space norm} & Euclidean & Frobenius\\
\hline {\bf sparsity inducing norm} & $\ell_1$ & nuclear\\
\hline {\bf dual norm} & $\ell_\infty$ & operator\\
\hline {\bf norm additivity} & disjoint support & orthogonal row and
column spaces\\
\hline {\bf convex optimization} & linear programming & semidefinite programming\\
\hline
\end{tabular}
\end{center}
\label{table:dictionary} \caption{\small A dictionary relating the
concepts of cardinality and rank minimization.}
\end{table}

\paragraph{Matrix vs.\  Vector Norms}

The three vector norms that play significant roles in the compressed
sensing framework are the $\ell_1$, $\ell_2$, and $\ell_\infty$ norms,
denoted by $\|x\|_1$, $\|x\|$ and $\|x\|_\infty$ respectively. These
norms have natural generalizations to matrices, inheriting many
appealing properties from the vector case.  In particular, there is a
parallel duality structure.

For a rectangular matrix $X \in \R^{m \times n}$, $\sigma_{i}(X)$
denotes the $i$-th largest singular value of $X$ and is equal to the
square-root of the $i$-th largest eigenvalue of $XX'$. The rank of
$X$ will usually be denoted by $r$, and is equal to the number of
nonzero singular values.  For matrices $X$ and $Y$ of the same
dimensions, we define the inner product in $\R^{m \times n}$ as
$\langle X,Y\rangle:=\trace(X'Y) = \sum_{i=1}^m \sum_{j=1}^n X_{ij}
Y_{ij}$. The norm associated with this inner product is called the
Frobenius (or Hilbert-Schmidt) norm $||\cdot||_F$. The Frobenius
norm is also equal to the Euclidean, or $\ell_2$, norm of the vector
of singular values, i.e.,
\[
    \|X\|_F := \sqrt{\langle X,X\rangle} =
    \sqrt{\trace(X'X)} = \left( \sum_{i=1}^m \sum_{j=1}^n X_{ij}^2 \right)^\frac{1}{2} =
    \left(\sum_{i=1}^{r} {\sigma_i^2}\right)^\frac{1}{2}.
\]
The operator norm (or induced 2-norm) of a matrix is equal to its
largest singular value (i.e., the $\ell_\infty$ norm of the singular
values):
\[
    \|X\| := \sigma_1(X).
\]
The nuclear norm of a matrix is equal to the sum of its singular
values, i.e.,
\[
\|X\|_* := \sum_{i=1}^{r} \sigma_i(X)\, ,
\]
and is alternatively known by several other names including the
Schatten $1$-norm, the Ky Fan $r$-norm, and the trace class norm.
Since the singular values are all positive, the nuclear norm is
equal to the $\ell_1$ norm of the vector of singular values. These
three norms are related by the following inequalities which hold for
any matrix $X$ of rank at most $r$:
\begin{equation}
||X|| \leq ||X||_F \leq ||X||_* \leq \sqrt{r} ||X||_F \leq r ||X||.
\label{eq:ineqnorms}
\end{equation}

\paragraph{Dual norms}
For any given norm $\|\cdot\|$ in an inner product space, there exists
a dual norm $\|\cdot\|_d$ defined as
\begin{equation}
    \|X\|_d := \max_Y \{ \langle X,Y \rangle \,:\, \|Y\|\leq 1\}.
\end{equation}
Furthermore, the norm dual to the norm $||\cdot||_d$ is again the
original norm $||\cdot||$.

In the case of vector norms in $\R^n$, it is well-known that the
dual norm of the $\ell_p$ norm (with $1 < p < \infty$) is the
$\ell_q$ norm, where $\frac{1}{p}+\frac{1}{q}=1$.  This fact is
essentially equivalent to H\"older's inequality.  Similarly, the
dual norm of the $\ell_\infty$ norm of a vector is the $\ell_1$
norm. These facts also extend to the matrix norms we have defined.
For instance, the dual norm of the Frobenius norm is the Frobenius
norm. This can be verified by simple calculus (or Cauchy-Schwarz),
since
\begin{equation*}
\max_{Y} \{\trace(X'Y) \,:\, \trace(Y'Y)\leq 1\}
\end{equation*}
is equal to $\|X\|_F$, with the maximizing $Y$ being equal to
$X/\|X\|_F$.  Similarly, as shown below, the dual norm of the operator
norm is the nuclear norm.  The proof of this fact will also allow us
to present variational characterizations of each of these norms as
semidefinite programs.

\begin{proposition}
The dual norm of the operator norm $||\cdot||$ in $\R^{m \times n}$
is the nuclear norm $||\cdot||_*$.
\label{prop:dualnorm}
\end{proposition}

\begin{proof}
First consider an $m \times n$ matrix $Z$.  The fact that $Z$ has
operator norm less than $t$ can be expressed as a linear matrix
inequality:
\begin{equation}
    \|Z\| \leq t
\quad \iff \quad t^2 \id_m - ZZ' \succeq 0 \quad \iff \quad
\left[\begin{array}{cc} t  \id_m & Z\\ Z' & t  \id_n
\end{array}\right] \succeq 0, \label{eq:charopnorm}
\end{equation}
where the last implication follows from a Schur complement argument.
As a consequence, we can give a semidefinite optimization
characterization of the operator norm, namely
\begin{equation}
\|Z\| = \min_t  \, t \qquad  \mbox{s.t.} \quad
\left[\begin{array}{cc} t \id_m & Z\\ Z' & t \id_n
    \end{array}\right] \succeq 0.
\end{equation}

Now let $X=U \Sigma V'$ be a singular value decomposition of an $m
\times n$ matrix $X$, where $U$ is an $m\times r$ matrix, $V$ is an
$n\times r$ matrix,  $\Sigma$ is a $r \times r$ diagonal matrix and
$r$ is the rank of $X$.  Let $Y := U \, V'$.  Then $\|Y\| = 1$ and
$\trace(XY') = \sum_{i=1}^{r} \sigma_i(X) = ||X||_*$, and hence the
dual norm is greater than or equal to the nuclear norm.

To provide an upper bound on the dual norm, we appeal to
semidefinite programming duality. From the characterization in
\eq{charopnorm}, the optimization problem
\[
    \max_Y \, \{ \langle X,Y \rangle \,:\, \|Y\|\leq 1\}
\]
is equivalent to the semidefinite program
\begin{equation}
\begin{array}{ll}
\displaystyle\max_Y & \trace(X'Y)\\
\mbox{s.t.} & \left[\begin{array}{cc} \id_m & Y\\ Y' & \id_n
    \end{array}\right] \succeq 0.
    \end{array}
\label{eq:nucnorm-sdp}
\end{equation}
The dual of this SDP (after an inconsequential rescaling) is given
by
\begin{equation}\label{eq:sdp-embedding}
\begin{array}{ll}
\displaystyle\min_{W_1,W_2} & \frac{1}{2}( \trace(W_1)+ \trace(W_2)) \\
\mbox{s.t.} & \left[\begin{array}{cc} W_1 & X\\ X' & W_2
    \end{array}\right] \succeq 0.
\end{array}
\end{equation}
Set $W_1:=U \Sigma U'$ and $W_2:=V \Sigma V'$.  Then the triple
$(W_1,W_2,X)$ is feasible for \eq{sdp-embedding} since
\[
\left[\begin{array}{cc} W_1 & X\\ X' & W_2
    \end{array}\right] = \left[\begin{array}{c} U \\ V \end{array}\right]
    \Sigma \left[\begin{array}{c} U
    \\ V \end{array}\right]' \succeq 0.
\]
Furthermore, we have $\trace(W_1) = \trace(W_2) = \trace(\Sigma)$, and
thus the objective function satisfies $(\trace(W_1)+\trace(W_2))/2 =
\trace \Sigma = \|X\|_*$. Since any feasible solution of
\eq{sdp-embedding} provides an upper bound for (\ref{eq:nucnorm-sdp}),
we have that the dual norm is less than or equal to the nuclear norm
of $X$, thus proving the proposition.
\end{proof}

Notice that the argument given in the proof above further shows that
the nuclear norm $||X||_*$ can be computed using either the SDP
\eq{nucnorm-sdp} or its dual \eq{sdp-embedding}, since there is no
duality gap between them. Alternatively, this could have also been
proven using a Slater-type interior point condition since both
\eq{nucnorm-sdp} and \eq{sdp-embedding} admit strictly feasible
solutions.

\paragraph{Convex envelopes of rank and cardinality functions}
Let $\mathcal{C}$ be a given convex set. The \emph{convex envelope}
of a (possibly nonconvex) function $f: \mathcal{C} \rightarrow \R$
is defined as the largest convex function $g$ such that $g(x)\leq
f(x)$ for all $x\in \mathcal{C}$ (see, e.g.,~\cite{HiriartBook}).
This means that among all convex functions, $g$ is the best
pointwise approximation to $f$.  In particular, if the optimal $g$
can be conveniently described, it can serve as an approximation to
$f$ that can be minimized efficiently.

By the chain of inequalities in~\eq{ineqnorms}, we have that
$\rank{X}\geq \|X\|_*/\|X\|$ for all $X$.  For all matrices with
$\|X\|\leq1$, we must have that $\rank{X}\geq \|X\|_*$, so the
nuclear norm is a convex lower bound of the rank function on the
unit ball in the operator norm.  In fact, it can be shown that this
is the tightest convex lower bound.
\begin{theorem}[\cite{FazelThesis}]\label{thm:convenv}
The convex envelope of $\rank{X}$ on the set $\{X\in\R^{m\times
n}~:~ \|X\|\leq 1\}$ is the nuclear norm $\|X\|_*$.
\end{theorem}
The proof is given in~\cite{FazelThesis} and uses a basic result
from convex analysis that establishes that (under some technical
conditions) the biconjugate of a function is its convex
envelope~\cite{HiriartBook}.

Theorem~\ref{thm:convenv} provides the following interpretation of
the nuclear norm heuristic for the affine rank minimization problem.
Suppose $X_0$ is the minimum rank solution of $\cA(X)=b$,
and $M=\|X_0\|$. The convex envelope of the rank on the set
$\mathcal{C}=\{X\in\R^{m\times n}~:~ \|X\|\leq M\}$ is $\|X\|_*/M$.
Let $X_*$ be the minimum nuclear norm solution of
$\cA(X)=b$. Then we have
\begin{equation*}
    \|X_*\|_*/M \leq \rank{X_0} \leq  \rank{X_*}
\end{equation*}
providing an upper and lower bound on the optimal rank when the norm
of the optimal solution is known.  Furthermore, this is the tightest
lower bound among all convex lower bounds of the rank function on
the set $\mathcal{C}$.

For vectors, we have a similar inequality.  Let $\card{x}$ denote
the cardinality function which counts the number of non-zero entries
in the vector $x$. Then we have $\card{x} \geq
\|x\|_1/\|x\|_\infty$. Not surprisingly, the $\ell_1$ norm is also
the convex envelope of the cardinality function over the set
$\{x\in\R^n~:~\|x\|_{\infty}\leq 1\}$. This result can be either
proven directly or can be seen as a special case of the above
theorem.

\paragraph{Additivity of rank and nuclear norm}
A function $f$ mapping a linear space $\mathcal{S}$ to $\R$ is
called \emph{subadditive} if $f(x+y)\leq f(x)+f(y)$. It is
\emph{additive} if $f(x+y) = f(x)+f(y)$.  In the case of vectors,
both the cardinality function and the $\ell_1$ norm are subadditive.
That is, if $x$ and $y$ are sparse vectors, then it always holds
that the number of non-zeros in $x+y$ is less than or equal to the
number of non-zeros in $x$ plus the number of non-zeros of $y$;
furthermore (by the triangle inequality) $\|x+y\|_1 \leq
\|x\|_1+\|y\|_1$.  In particular, the cardinality function is
additive exactly when the vectors $x$ and $y$ have disjoint support.
In this case, the $\ell_1$ norm is also additive, in the sense that
$\|x+y\|_1 = \|x\|_1+\|y\|_1$.

For matrices, the rank function is subadditive.  For the rank to be
additive, it is necessary and sufficient that the row and column
spaces of the two matrices intersect only at the origin, since in
this case they operate in essentially disjoint spaces (see, e.g.,
\cite{MarsagliaStyan}).  As we will show below, a related condition
that ensures that the nuclear norm is additive, is that the matrices
$A$ and $B$ have row and column spaces that are \emph{orthogonal}.
In fact, a compact sufficient condition for the additivity of the
nuclear norm will be that $AB'=0$ and $A'B=0$.  This is a stronger
requirement than the aforementioned condition for rank additivity,
as orthogonal subspaces only intersect at the origin. The disparity
arises because the nuclear norm of a linear map depends on the
choice of the inner products on the spaces $\R^m$ and $\R^n$ on
which the matrix acts, whereas the rank is independent of such a
choice.

\begin{lemma}\label{lemma:l1-split}
Let $A$ and $B$ be matrices of the same dimensions.  If $AB'=0$ and
$A'B=0$ then $\|A+B\|_* = \|A\|_* + \|B\|_*$.
\end{lemma}

\begin{proof}
Partition the singular value decompositions of $A$ and $B$ to
reflect the zero and non-zero singular vectors
\begin{equation*}
    A=\left[\begin{array}{cc}U_{A1} & U_{A2}\end{array}\right]
   \left[\begin{array}{cc}\Sigma_A & \\ & 0\end{array}\right]
   \left[\begin{array}{cc}V_{A1} & V_{A2}\end{array}\right]'
   \qquad
    B=\left[\begin{array}{cc}U_{B1} & U_{B2}\end{array}\right]
   \left[\begin{array}{cc}\Sigma_B & \\ & 0\end{array}\right]
   \left[\begin{array}{cc}V_{B1} & V_{B2}\end{array}\right]'\,.
\end{equation*}
The condition $AB'=0$ implies that $V_{A1}'V_{B1}=0$, and similarly,
$A'B=0$ implies that $U_{A1}'U_{B1}=0$.  Hence, there exist matrices
$U_C$ and $V_C$ such that $[U_{A1}\, U_{B1}\, U_C]$ and $[V_{A1}\,
V_{B1}\, V_C]$ are orthogonal matrices. Thus, the following are
valid singular value decompositions for $A$ and $B$:
\begin{align*}
    A&=\left[\begin{array}{ccc}U_{A1} & U_{B1} & U_C\end{array}\right]
   \left[\begin{array}{ccc}\Sigma_A & & \\ & 0 & \\ & & 0\end{array}\right]
   \left[\begin{array}{ccc}V_{A1} & V_{B1} & V_C\end{array}\right]'
   \\
    B&=\left[\begin{array}{ccc}U_{A1} & U_{B1} & U_C\end{array}\right]
   \left[\begin{array}{ccc}0 & & \\ & \Sigma_B & \\ & & 0\end{array}\right]
   \left[\begin{array}{ccc}V_{A1} & V_{B1} & V_C\end{array}\right]'\,.
\end{align*}
In particular, we have that
\begin{equation*}
    A+B = \left[\begin{array}{cc}U_{A1} & U_{B1}\end{array}\right]
   \left[\begin{array}{cc}\Sigma_A & \\ & \Sigma_B\end{array}\right]
   \left[\begin{array}{cc}V_{A1} & V_{B1}\end{array}\right]'\,.
\end{equation*}
This shows that the singular values of $A+B$ are equal to the union
(with repetition) of the singular values of $A$ and $B$. Hence,
$\|A+B\|_*=\|A\|_*+\|B\|_*$ as desired.
\end{proof}

\begin{corollary}
Let $A$ and $B$ be matrices of the same dimensions.  If the row and
column spaces of $A$ and $B$ are orthogonal, then $\|A+B\|_* =
\|A\|_*+\|B\|_*$.
\end{corollary}
\begin{proof}
It suffices to show that if the row and column spaces of $A$ and $B$
are orthogonal, then $AB'=0$ and $A'B=0$. But this is immediate: if
the columns of $A$ are orthogonal to the columns of $B$, we have
$A'B=0$. Similarly, orthogonal row spaces imply that $AB'=0$ as
well.
\end{proof}

\paragraph{Nuclear norm minimization}
Let us turn now to the study of equality-constrained norm
minimization problems where we are searching for a matrix $X \in
\R^{m \times n}$ of minimum nuclear norm belonging to a given affine
subspace. In our applications, the subspace is usually described by
linear equations of the form $\cA(X)= b$, where
$\cA:\R^{m \times n} \rightarrow \R^p$ is a linear mapping.
This problem admits the primal-dual convex formulation
\begin{equation}
    \begin{aligned}
    \min_{X}&\quad \|X\|_*          & \qquad \qquad \qquad \qquad    \max_z & \quad b'z \\
    \mbox{s.t} &\quad \cA(X) = b &     \mbox{s.t.} & \quad ||\cA^*(z)|| \leq 1,
    \end{aligned}
\label{eq:normminimization}
\end{equation}
where $\cA^*:\R^p \rightarrow \R^{m \times n}$ is the
adjoint of $\cA$.  The formulation
(\ref{eq:normminimization}) is valid for any norm minimization
problem, by replacing the norms appearing above by any dual pair of
norms.  In particular,  if we replace the nuclear norm with the
$\ell_1$ norm and the operator norm with the $\ell_\infty$ norm, we
obtain a primal-dual pair of optimization problems, that can be
reformulated in terms of linear programming.

Using the SDP characterizations of the nuclear and operator norms
given in \eq{nucnorm-sdp}-\eq{sdp-embedding} above allows us to
rewrite \eq{normminimization} as the primal-dual pair of
semidefinite programs
\begin{equation}
\label{eq:nucnorm-sdp-pd}
\begin{aligned}
\min_{X,W_1,W_2} && \textstyle{\frac{1}{2}}(\trace(W_1)+\trace&(W_2)) &
\qquad \qquad \qquad \max_{z} & \quad b'z
\\
\mbox{s.t.} && \left[\begin{array}{cc} W_1 & X\\ X' & W_2
    \end{array}\right] &\succeq 0  &
\mbox{s.t.} & \quad \left[\begin{array}{cc} \id_m & \cA^* (z)\\
\cA^*(z)' & \id_n
    \end{array}\right] \succeq 0. \\
&& \cA(X) &= b
\end{aligned}
\end{equation}

\paragraph{Optimality conditions}
In order to describe the optimality conditions for the norm
minimization problem \eq{normminimization}, we must first
characterize the subdifferential of the nuclear norm. Recall that
for a convex function $f:\R^{n} \rightarrow \R$, the subdifferential
of $f$ at $x \in \R^n$ is the compact convex set
\[
\partial f(x) := \{ d \in \R^n \, : \, f(y) \geq f(x) + \langle d , y-x \rangle
\quad \forall y \in \R^n \}.
\]
Let $X$ be an $m\times n$ matrix with rank $r$ and let $X=U\Sigma
V'$ be a singular value decomposition where $U \in \R^{m\times r}$,
$V \in \R^{n\times r}$ and $\Sigma$ is an $r\times r$ diagonal
matrix.  The subdifferential of the nuclear norm at $X$ is then
given by (see, e.g., \cite{Watson92})
\begin{equation}
    \partial \|X\|_* = \{UV' + W \, : \, W \mbox{ and } X
\mbox{ have orthogonal row and column spaces, and } \|W\|\leq 1 \}.
\label{eq:subdiffnuc}
\end{equation}
For comparison, recall the case of the $\ell_1$ norm, where $T$
denotes the support of the $n$-vector $x$, $T^c$ is the complement of
$T$ in the set $\{1,\ldots,n\}$, and
\begin{equation}
    \partial \|x\|_1 = \{d \in \R^n \, : \, d_i = \sign(x) \mbox{ for } i\in
    T, \quad |d_i|\leq 1 \mbox{ for } i \in T^c\}.
\label{eq:subdiffl1}
\end{equation}
The similarity between \eq{subdiffnuc} and \eq{subdiffl1} is
particularly transparent if we recall the \emph{polar decomposition}
of a matrix into a product of orthogonal and positive semidefinite
matrices  (see, e.g., \cite{HornJohnsonBook2}).  The ``angular''
component of the matrix $X$ is exactly given by $UV'$. Thus, these
subgradients always have the form of an ``angle'' (or sign), plus
possibly a contraction in an orthogonal direction if the norm is not
differentiable at the current point.

We can now write concise optimality conditions for the optimization
problem \eq{normminimization}. A matrix $X$ is optimal for
\eq{normminimization} if there exists a vector $z\in\R^p$ such that
\begin{equation}
    \cA(X) = b, \qquad \qquad
    \cA^*(z)  \in \partial \|X\|_*.
\label{eq:optsubdiff}
\end{equation}
The first condition in \eq{optsubdiff} requires feasibility of the
linear equations, and the second one guarantees that there is no
feasible direction of improvement. Indeed, since $\cA^*(z)$ is
in the subdifferential at $X$, for any $Y$ in the primal feasible set
of \eq{normminimization} we have
\[
||Y||_* \geq ||X||_* + \langle \cA^*(z), Y-X \rangle =
||X||_* + \langle z , \cA(Y-X) \rangle =
||X||_*,
\]
where the last step follows from the feasibility of $X$ and $Y$.  As
we can see, the optimality conditions \eq{optsubdiff} for the nuclear
norm minimization problem exactly parallel those of the $\ell_1$
optimization case.

These optimality conditions can be used to check and certify whether
a given candidate $X$ is indeed the minimum nuclear norm solution.
For this, it is sufficient (and necessary) to find a vector $z \in
\R^p$ in the subdifferential of the norm, i.e., such that the left-
and right-singular spaces of $\cA^*(z)$ are aligned with
those of $X$, and is a contraction in the orthogonal complement.

\section{Restricted Isometry and Recovery of Low-Rank Matrices}
\label{sec:rip}

Let us now turn to the central problem analyzed in this paper. Let
$\cA:\R^{m \times n} \rightarrow \R^p$ be a linear map and
let $X_0$ be a matrix of rank $r$.  Set $b:=\cA(X_0)$, and
define the convex optimization problem
\begin{equation}\label{eq:nucnorm-prob}
      X^* :=  \arg\min_{X} \|X\|_*\qquad \mbox{s.t.}\quad \cA(X) = b.
\end{equation}
In this section, we will characterize specific cases when we can
\emph{a priori} guarantee that $X^*=X_0$. The key conditions will be
determined by the values of a sequence of parameters $\delta_r$ that
quantify the behavior of the linear map $\cA$ when restricted
to the subvariety of matrices of rank $r$. The following definition is
the natural generalization of the Restricted Isometry Property from
vectors to matrices.

\begin{definition} Let $\cA:\R^{m \times n} \rightarrow \R^p$ be a linear map.
Without loss of generality, assume $m\leq n$.  For every integer $r$
with $1\leq r \leq m$, define the $r$-restricted isometry constant
to be the smallest number $\delta_r(\cA)$ such that
\begin{equation}\label{eq:RIP}
    (1-\delta_r(\cA))\|X\|_F \leq \|\cA(X)\| \leq
    (1+\delta_r(\cA))\|X\|_F
\end{equation}
holds for all matrices $X$ of rank at most $r$.
\end{definition}

Note that by definition, $\delta_r(\cA) \leq
\delta_{r^\prime}(\cA)$ for $r \leq r^\prime$.

The Restricted Isometry Property for sparse vectors was developed by
Cand\`es and Tao in~\cite{Candes05}, and requires \eq{RIP} to hold
with the Euclidean norm replacing the Frobenius norm and rank being
replaced by cardinality. Since for diagonal matrices, the Frobenius
norm is equal to the Euclidean norm of the diagonal, this definition
reduces to the original Restricted Isometry Property
of~\cite{Candes05} in the diagonal
case.\footnote{In~\cite{Candes05}, the authors define the restricted
isometry properties with squared norms.  We note here that the
analysis is identical modulo some algebraic rescaling of constants.
We choose to drop the squares as it greatly simplifies the analysis
in Section~\ref{sec:gaussian}.}

Unlike the case of ``standard'' compressed sensing, our RIP
condition for low-rank matrices cannot be interpreted as
guaranteeing all sub-matrices of the linear transform $\cA$
of a certain size are well conditioned.  Indeed, the set of matrices
$X$ for which~\eq{RIP} must hold is \emph{not} a finite union of
subspaces, but rather a certain ``generalized Stiefel manifold,''
which is also an algebraic variety (in fact, it is the $r$th-secant
variety of the variety of rank-one matrices). Surprisingly, we are
still able to derive analogous recovery results for low-rank
solutions of equations when $\cA$ obeys this RIP condition.
Furthermore, we will see in Section~\ref{sec:gaussian} that many
ensembles of random matrices have the Restricted Isometry Property
with $\delta_r$ quite small with high probability for reasonable
values of $m$,$n$, and $p$.

The following two recovery theorems will characterize the power of
the restricted isometry constants.  Both theorems are more or less
immediate generalizations from the sparse case to the low-rank case
and use only minimal properties of the rank of matrices and the
nuclear norm. The first theorem generalizes Lemma 1.3
in~\cite{Candes05} to low-rank recovery.
\begin{theorem}\label{thm:reclowrank}
Suppose that $\delta_{2r}<1$ for some integer $r\geq 1$.  Then $X_0$
is the only matrix of rank at most $r$ satisfying
$\cA(X)=b$.
\end{theorem}

\begin{proof}
Assume, on the contrary, that there exists a rank $r$ matrix $X$
satisfying $\cA(X)=b$ and $X \neq X_0$.  Then $Z:=X_0-X$ is
a nonzero matrix of rank at most $2r$, and $\cA(Z)=0$. But
then we would have $0=\|\cA(Z)\| \geq (1-\delta_{2r})\|Z\|_F
> 0$ which is a contradiction.
\end{proof}

The proof of the preceding theorem is identical to the argument
given by Cand\`es and Tao and is an immediate consequence of our
definition of the constant $\delta_r$. No adjustment is necessary in
the transition from sparse vectors to low-rank matrices.  The key
property used is the sub-additivity of the rank.

Next, we state a weak $\ell_1$-type recovery theorem whose proof
mimics the approach in~\cite{Candes05b}, but for which a few details
need to be adjusted when switching from vectors to matrices.

\begin{theorem}\label{thm:recweak}
Suppose that $r\geq 1$ is such that $\delta_{5r} < 1/10$. Then
$X^*=X_0$.
\end{theorem}

We will need the following technical lemma that shows for any two
matrices $A$ and $B$, we can decompose $B$ as the sum of two matrices
$B_1$ and $B_2$ such that $\rank{B_1}$ is not too large and $B_2$
satisfies the conditions of Lemma~\ref{lemma:l1-split}.  This will be
the key decomposition for proving Theorem~\ref{thm:recweak}.

\begin{lemma}\label{lemma:rank-partition}
Let $A$ and $B$ be matrices of the same dimensions.  Then there
exist matrices $B_1$ and $B_2$ such that
\begin{enumerate}
\item   $B = B_1+B_2$
\item \label{rankcon}  $\rank{B_1} \leq 2 \rank{A}$
\item   $A B_2'=0$ and $A 'B_2=0$
\item   $\langle B_1, B_2 \rangle=0$
\end{enumerate}
\end{lemma}
\begin{proof}
Consider a full singular value decomposition of $A$
\begin{equation*}
    A = U\left[\begin{array}{cc} \Sigma & 0 \\ 0 &
    0\end{array}\right] V',
\end{equation*}
and let $\hat{B} := U'BV$.  Partition $\hat{B}$ as
\begin{equation*}
    \hat{B} = \left[\begin{array}{cc} \hat{B}_{11} & \hat{B}_{12}\\ \hat{B}_{21} &  \hat{B}_{22}  \end{array}  \right].
\end{equation*}
Defining now
\begin{equation*}
    {B}_1 := U \left[\begin{array}{cc} \hat{B}_{11} & \hat{B}_{12}\\ \hat{B}_{21} &  0  \end{array}
    \right] V', \qquad \qquad
    {B}_2 := U \left[\begin{array}{cc} 0 & 0\\ 0 & \hat{B}_{22} \end{array}  \right] V',
\end{equation*}
it can be easily verified that $B_1$ and $B_2$ satisfy the
conditions (1)--(4).
\end{proof}

We now proceed to a proof of Theorem~\ref{thm:recweak}.

\begin{proof}{\bf [of Theorem~\ref{thm:recweak}]}
By optimality of $X^*$, we have $\|X_0\|_*\geq\|X^*\|_*$. Let
$R:=X^*-X_0$. Applying Lemma~\ref{lemma:rank-partition} to the
matrices $X_0$ and $R$, there exist matrices $R_0$ and $R_c$ such that
$R=R_0+R_c$, $\rank{R_0} \leq 2 \rank{X_0}$, and $X_0 R_c'=0$ and
$X_0' R_c=0$.  Then,
\begin{equation}
\begin{split}
    \|X_0\|_* \geq \|X_0+R\|_* \geq \|X_0+R_c\|_* - \|R_0\|_* = \|X_0\|_* + \|R_c\|_* -
    \|R_0\|_*,
\end{split}
\end{equation}
where the middle assertion follows from the triangle inequality and
the last one from Lemma~\ref{lemma:l1-split}. Rearranging terms, we
can conclude that
\begin{equation}
    \|R_0\|_* \geq \|R_c\|_*.
\label{eq:r0rc}
\end{equation}

Next we partition $R_c$ into a sum of matrices $R_1, R_2, \ldots$,
each of rank at most $3r$.  Let $R_c= U\diag(\sigma)V'$ be the
singular value decomposition of $R_c$.  For each $i\geq 1$ define
the index set $I_i = \{3r(i-1)+1,\ldots,3ri\}$, and let
$R_i:=U_{I_i}\diag(\sigma_{I_i})V_{I_i}'$ (notice that $\langle R_i,
R_j \rangle = 0$ if $i \not = j$). By construction, we have
\begin{equation}
    \sigma_k \leq \frac{1}{3r} \sum_{j\in I_i} \sigma_j
    \qquad \forall\, k\in I_{i+1},
\end{equation}
which implies $\|R_{i+1}\|_F^2 \leq \frac{1}{3r} \|R_i\|_*^2$.  We
can then compute the following bound
\begin{equation}
    \sum_{j\geq 2} \|R_j\|_F \leq \frac{1}{\sqrt{3r}} \sum_{j\geq 1}
    \|R_j\|_* = \frac{1}{\sqrt{3r}} \|R_c\|_* \leq \frac{1}{\sqrt{3r}}
    \|R_0\|_* \leq \frac{\sqrt{2r}}{\sqrt{3r}}\|R_0\|_F\,,
\end{equation}
where the last inequality follows from~\eq{ineqnorms} and the fact
that $\rank{R_0} \leq 2r$.  Finally, note that the rank of $R_0+R_1$
is at most $5r$, so we may put this all together as
\begin{equation}\label{eq:l1rec}
\begin{split}
    \|\cA(R)\| &\geq \|\cA(R_0+R_1)\| -
    \sum_{j\geq 2} \|\cA(R_j)\| \\
    & \geq (1-\delta_{5r}) \, \|R_0+R_1\|_F -
    (1+\delta_{3r}) \, \sum_{j\geq 2} \|R_j\|_F\\
    & \geq \left((1-\delta_{5r}) -  \sqrt{\tfrac{2}{3}}(1+\delta_{3r}) \right)
    \|R_0\|_F\\
    & \geq \left((1-\delta_{5r}) -  \tfrac{9}{11}(1+\delta_{3r}) \right)
    \|R_0\|_F.
\end{split}
\end{equation}
By assumption $\cA(R) = \cA(X^*-X_0) = 0$, so if the
factor on the right-hand side is strictly positive, $R_0=0$, which
further implies $R_c=0$ by \eq{r0rc}, and thus $X^*=X_0$. Simple
algebra reveals that the right-hand side is positive when $9
\delta_{3r}+ 11 \delta_{5r} < 2$. Since $\delta_{3r} \leq
\delta_{5r}$, we immediately have that $X^* = X_0$ if $\delta_{5r} <
1/10$.
\end{proof}

The rational number ($9/11$) in the proof of the theorem is chosen
for notational simplicity and is clearly not optimal. A slightly
tighter bound can be achieved working directly with the second to
last line in \eq{l1rec}.  The most important point is that our
recovery condition on $\delta_{5r}$ is an absolute constant,
independent of $m$, $n$, $r$, and $p$.

We have yet to demonstrate any specific linear mappings
$\cA$ for which $\delta_r<1$. We shall show in the next
section that linear transformations sampled from several families of
random matrices with appropriately chosen dimensions have this
property with overwhelming probability. The analysis is again
similar to the compressive sampling literature, but several details
specific to the rank recovery problem need to be employed.

\section{Nearly Isometric Families}
\label{sec:gaussian}

In this section, we will demonstrate that when we sample linear maps
from a class of probability distributions obeying certain tail
bounds, then they will obey the Restricted Isometry Property
\eq{RIP} as $p$, $m$, and $n$ tend to infinity at appropriate rates.
The following definition characterizes this family of random linear
transformation.

\begin{definition}
Let $\cA$ be a random variable that takes values in linear maps from
$\R^{m \times n}$ to $\R^p$.  We say that $\cA$ is \emph{nearly
isometrically distributed} if for all $X \in \R^{m\times n}$
\begin{equation}
    \mathbf{E}[ \|\cA(X)\|^2] = \|X\|_F^2\,
\end{equation}
and for all $0<\epsilon<1$ we have,
\begin{equation}\label{eq:concentration-point}
    \mathbf{P}( |\|\cA(X)\|^2 - \|X\|_F^2| \geq \epsilon\|X\|_F^2 ) \leq
        2 \exp\left(-\frac{p}{2}({\epsilon^2}/{2}-{\epsilon^3}/{3})\right)
\end{equation}
and for all $t>0$, we have
\begin{equation}\label{eq:concentration-norm}
    \mathbf{P}\left( \|\cA\| \geq 1 + \sqrt{\frac{mn}{p}} + t \right) \leq \exp(-\gamma p t^2)
\end{equation}
for some constant $\gamma>0$.
\end{definition}

There are two ingredients for a random linear map to be nearly
isometric. First, it must be isometric in expectation.  Second, the
probability of large distortions of length must be exponentially
small. The exponential bound in \eq{concentration-point} guarantees
union bounds will be small even for rather large sets. This
concentration is the typical ingredient required to prove the
Johnson-Lindenstrauss Lemma (cf~\cite{Achlioptas03,Dasgupta03}).

The majority of nearly isometric random maps are described in terms
of random matrices.  For a linear map $\cA:\R^{m\times
n}\rightarrow\R^p$, we can always write its matrix representation as
\begin{equation}\label{eq:gauss-lin}
    \cA(X) = \mathbf{A}\vvec{X}\, ,
\end{equation}
where $\vvec{X}$ denotes the vector of $X$ with its columns stacked
in order on top of one another, and $\mathbf{A}$ is a $p \times mn$
matrix.  We now give several examples of nearly isometric random
variables in this matrix representation. The most well known is the
ensemble with independent, identically distributed (i.i.d.) Gaussian
entries~\cite{Dasgupta03}
\begin{equation}
    A_{ij} \sim \mathcal{N}(0,\frac{1}{p})\,.
\end{equation}
We also mention the two following ensembles of matrices, described
in~\cite{Achlioptas03}. One has entries sampled from an i.i.d.
symmetric Bernoulli distribution
\begin{equation}
        A_{ij} = \begin{cases}\sqrt{\frac{1}{p}} & \mbox{with probability } \frac{1}{2}\\
                                -\sqrt{\frac{1}{p}} & \mbox{with probability } \frac{1}{2}
        \end{cases}\, ,
\end{equation}
and the other has zeros in two-thirds of the entries
\begin{equation}
        A_{ij} = \begin{cases}
\sqrt{\frac{3}{p}} & \mbox{with probability } \frac{1}{6}\\
0 & \mbox{with probability } \frac{2}{3}\\
-\sqrt{\frac{3}{p}} & \mbox{with probability } \frac{1}{6}
        \end{cases}\, .
\end{equation}
The fact that the top singular value of the matrix $\mathbf{A}$ is
concentrated around $1+\sqrt{D/p}$ for all of these ensembles
follows from the work of Yin, Bai, and Krishnaiah, who showed that
whenever the entries $A_{ij}$ are i.i.d. with zero mean and finite
fourth moment, then the maximum singular value of $\mathbf{A}$ is
almost surely $1+\sqrt{D/p}$ for $D$ sufficiently
large~\cite{Bai88}. El Karoui uses this result to prove the
concentration inequality \eq{concentration-norm} for all such
distributions~\cite{ElKarouiThesis}.  The result for Gaussians is
rather tight with $\gamma=1/2$ (see, e.g., \cite{Davidson01}).

Finally, note that a random projection also obeys all of the necessary
concentration inequalities. Indeed, since the norm of a random
projection is exactly $\sqrt{D/p}$, \eq{concentration-norm} holds
trivially. The concentration inequality \eq{concentration-point} is
proven in~\cite{Dasgupta03}.

The main result of this section is the following:
\begin{theorem}\label{thm:ripscaling}
Fix $0 < \delta < 1$.  If $\cA$ is a nearly isometric random
variable, then for every $1\leq r \leq m$, there exist constants
$c_0,c_1>0$ depending only on $\delta$ such that, with probability
at least $1-\exp(-c_1 p)$, $\delta_r(\cA) \leq \delta$
whenever $p\geq c_0 r(m+n)\log(mn)$.
\end{theorem}

The proof will make use of standard techniques in concentration of
measure. We first extend the concentration results of~\cite{Wakin07}
to subspaces of matrices.  We will show that the distortion of a
subspace by a linear map is robust to perturbations of the
subspace. Finally, we will provide an epsilon net over the set of all
subspaces and, using a union bound, will show that with overwhelming
probability, nearly isometric random variables will obey the
Restricted Isometry Property \eq{RIP} as the size of the matrices tend
to infinity.

The following lemma characterizes the behavior of a nearly isometric
random mapping $\cA$ when restricted to an arbitrary subspace
of matrices $U$ of dimension $d$.
\begin{lemma}\label{lemma:lipshitz}
Let $\cA$ be a nearly isometric linear map and let $U$ be an
arbitrary subspace of $m\times n$ matrices with $d = \dim(U) \leq p$.
Then for any $0 < \delta < 1$ we have
\begin{equation}\label{eq:ssrip}
    (1-\delta) \|X\|_F \leq \|\cA(X)\| \leq
    (1+\delta)\|X\|_F \qquad \forall\, X\in U
\end{equation}
with probability at least
\begin{equation}\label{eq:ssrip-prob}
    1 - 2(12/\delta)^d\exp\left(-\frac{p}{2}(\delta^2/8 -
    \delta^3/24)\right)\,.
\end{equation}
\end{lemma}

\begin{proof}
The proof of this theorem is identical to the argument
in~\cite{Wakin07} where the authors restricted their attention to
subspaces aligned with the coordinate axes. We will sketch the proof
here as the argument is straightforward.

There exists a finite set $\Omega$ of at most $(12/\delta)^d$ points
such that for every $X\in U$ with $\|X\|_F \leq 1$, there exists a $Q
\in \Omega$ such that $\|X-Q\|_F \leq \delta/4$.  By the standard
union bound, the concentration inequality \eq{concentration-point}
holds for all $Q\in\Omega$ with $\epsilon = \delta/2$ with probability
at least \eq{ssrip-prob}.  If \eq{concentration-point} holds for all
$Q\in \Omega$, then we immediately have that $(1-\delta/2)\|Q\|_F\leq
\|\cA(Q)\|\leq (1+\delta/2)\|Q\|_F$ for all $Q\in\Omega$ as
well.

Let $X$ be in $\{X \in U: \|X\|_F \leq 1\}$, and $M$ be the maximum
of $\|\cA(X)\|$ on this set. Then there exists a $Q\in
\Omega$ such that $\|X-Q\|_F\leq \delta/4$. We then have
\[
    \|\cA(X)\| \leq \|\cA(Q)\| +
    \|\cA(X-Q)\| \leq 1+\delta/2+ M \delta/4,
\]
and since $M \leq 1+\delta/2+M\delta/4$ by definition, we have
$M \leq 1+\delta$.  The lower bound is proven by the following chain of
inequalities
\[
    \|\cA(X)\|  \geq \|\cA(Q)\| -
    \|\cA(X-Q)\| \geq 1 - \delta/2 - (1+\delta)\delta/4 \geq
    1-\delta.
\]
\end{proof}

The proof of preceding lemma revealed that the near isometry of a
linear map is robust to small perturbations of the matrix on which
the map is acting.  We will now show that this behavior is robust
with respect to small perturbations of the subspace $U$ as well.
This perturbation will be measured in the natural distance between
two subspaces
\begin{equation}
    \rho(T_1,T_2) := \|P_{T_1}-P_{T_2}\|,
\label{eq:metricG}
\end{equation}
where $T_1$ and $T_2$ are subspaces and $P_{T_i}$ is the orthogonal
projection associated with each subspace. This distance measures the
operator norm of the difference between the corresponding projections,
and is equal to the sine of the largest principal angle between $T_1$
and $T_2$ \cite{EdelmanAngle}.

The set of all $d$-dimensional subspaces of $\R^D$ is commonly known
as the Grassmannian manifold $\mathfrak{G}(D,d)$.  We will endow it
with the metric $\rho(\cdot,\cdot)$ given by~(\ref{eq:metricG}), also
known as the \emph{projection 2-norm}. In the following lemma we
characterize and quantify the change in the isometry constant $\delta$
as one smoothly moves through the Grassmannian.

\begin{lemma}
Let $U_1$ and $U_2$ be $d$-dimensional subspaces of $\R^D$.
Suppose that for all $X\in U_1$,
\begin{equation}
    (1-\delta)\|X\|_F \leq \|\cA(X)\| \leq
    (1+\delta)\|X\|_F
\end{equation}
for some constant $0<\delta<1$.  Then for all $Y \in U_2$
\begin{equation}
    (1-\delta')\|Y\|_F \leq \|\cA(Y)\| \leq
    (1+\delta')\|Y\|_F
\end{equation}
with
\begin{equation}
\delta' = \delta + (1+\|\cA\|) \cdot \rho(U_1,U_2)\,.
\end{equation}
\end{lemma}

\begin{proof}
Consider any $Y \in U_2$. Then
\begin{equation}
\begin{split}
    \|\cA(Y)\| &=
    \left\|\cA\left(P_{U_1}(Y) -
    [P_{U_1}-P_{U_2}](Y)\right)\right\|\\
    &\leq
    \left\|\cA(P_{U_1}(Y))\| +
    \|\cA\left([P_{U_1}-P_{U_2}](Y)\right)\right\|\\
    &\leq (1+\delta) \|P_{U_1}(Y)\|_F
    + \|\cA\| \|P_{U_1}-P_{U_2}\|\|Y\|_F\\
    &\leq \left(1+\delta + \|\cA\| \|P_{U_1}-P_{U_2}\| \right)\|Y\|_F.
\end{split}
\end{equation}

Similarly, we have
\begin{equation}
\begin{split}
    \|\cA(Y)\| &\geq
    \left\|\cA(P_{U_1}(Y))\| -
    \|\cA\left([P_{U_1}-P_{U_2}](Y)\right)\right\|\\
    &\geq (1-\delta) \|P_{U_1}(Y)\|_F
    - \|\cA\|\|P_{U_1}-P_{U_2}\|\|Y\|_F\\
    &\geq (1-\delta) \|Y\|_F - (1-\delta)\|(P_{U_1}-P_{U_2})(Y)\|_F
    - \|\cA\| \|P_{U_1}-P_{U_2}\| \|Y\|_F\\
    &\geq \left[1- \delta  - (\|\cA\| +
    1)\|P_{U_1}-P_{U_2}\|\right]\|Y\|_F,
\end{split}
\end{equation}
which completes the proof.
\end{proof}

To apply these concentration results to low-rank matrices, we
characterize the set of all matrices of rank at most $r$ as a union
of subspaces.  Let $V \subset \R^m$ and $W \subset \R^n$ be fixed
subspaces of dimension $r$. Then the set of all $m\times n$ matrices
$X$ whose row space is contained in $W$ and column space is
contained in $V$ forms an $r^2$-dimensional subspace of matrices of
rank less than or equal to $r$. Denote this subspace as $\Sigma(V,W)
\subset \R^{m \times n}$. Any matrix of rank less than or equal to
$r$ is an element of some $\Sigma(V,W)$ for a suitable pair of
subspaces, i.e., the set
\[
\Sigma_{mnr} := \{ \Sigma(V,W) \quad : \quad V \in \mathfrak{G}(m,r),
\quad W \in \ \mathfrak{G}(n,r) \}.
\]
We now characterize how many subspaces are necessary to cover this
set to arbitrary resolution. The \emph{covering number}
$\mathfrak{N}(\epsilon)$ of $\Sigma_{mnr}$ at resolution $\epsilon$
is defined to be the smallest number of subspaces $(V_i,W_i)$ such
that for any pair of subspaces $(V,W)$, there is an $i$ with
$\rho(\Sigma(V,W),\Sigma(V_i,W_i))\leq \epsilon$. That is, the
covering number is the smallest cardinality of an $\epsilon$-net.
The following Lemma characterizes the cardinality of such a set.

\begin{lemma}\label{lemma:covering}
The covering number $\mathfrak{N}(\epsilon)$ of $\Sigma_{mnr}$ is
bounded above by
\begin{equation}
\mathfrak{N}(\epsilon) \leq
\left(\frac{{2}C_0}{\epsilon}\right)^{r(m+n-2r)}
\end{equation}
where  $C_0$ is a constant independent of $\epsilon$, $m$, $n$, and
$r$.
\end{lemma}
\begin{proof}
Note that the projection operator onto $\Sigma(V,W)$ can be written
as $P_{\Sigma(V,W)} = P_{V} \otimes P_{W}$, so for a pair of
subspaces $(V_1,W_1)$ and $(V_2,W_2)$, we have
\begin{equation}
\begin{split}
\rho(\Sigma(V_1,W_1),\Sigma(V_2,W_2)) &= \| P_{\Sigma(V_1,W_1)}-
P_{\Sigma(V_2,W_2)}\| \\
&=\|P_{V_1}\otimes P_{W_1}- P_{V_2}\otimes P_{W_2}\| \\
&=\| (P_{V_1}-P_{V_2})\otimes P_{W_1} + P_{V_2}\otimes
(P_{W_1} - P_{W_2})\| \\
&\leq \|P_{V_1}-P_{V_2}\| \|P_{W_1}\| +
\|P_{V_2}\| \|P_{W_1} - P_{W_2}\| \\
&=\rho(V_1,V_2) +\rho(W_1,W_2).
\end{split}
\end{equation}

The conditions $\rho(V_1,V_2)\leq \frac{\epsilon}{{2}}$ and
$\rho(W_1, W_2)\leq\frac{\epsilon}{{2}}$ together imply that
$\rho(\Sigma(V_1,W_1),\Sigma(V_2,W_2)) \leq \rho(V_1,V_2)+ \rho(W_1,
W_2)\leq \epsilon$. Let $V_1,\ldots, V_{N_1}$ cover the set of
$r$-dimensional subspaces of $\R^m$ to resolution $\epsilon/{2}$ and
$U_1, \ldots, U_{N_2}$ cover the $r$-dimensional subspaces of $\R^n$
to resolution $\epsilon/{2}$. Then for any $(V,W)$, there exist $i$
and $j$ such that $\rho(V,V_i)\leq \epsilon/2$ and $\rho(W,W_j)\leq
\epsilon/2$. Therefore, $\mathfrak{N}(\epsilon)\leq N_1 N_2$. By the
work of Szarek on $\epsilon$-nets of the Grassmannian
(\cite{Szarek83}, \cite[Th. 8]{SzarekHomog}) there is a universal
constant $C_0$, independent of $m$, $n$, and $r$, such that
\begin{equation}
    N_1 \leq \left(\frac{{2}C_0}{\epsilon}\right)^{r(m-r)}\qquad
    \mbox{and} \qquad
    N_2 \leq \left(\frac{{2}C_0}{\epsilon}\right)^{r(n-r)}
\end{equation}
which completes the proof.
\end{proof}

The exact value of the universal constant $C_0$ is not provided by
Szarek in \cite{SzarekHomog}.  It takes the same value for any
homogeneous space whose automorphism group is a subgroup of the
orthogonal group, and is independent of the dimension of the
homogeneous space. Hence, one might expect this constant to be quite
large.  However, it is known that for the sphere $C_0\leq
3$~\cite{LorentzBook}, and there is no indication that this constant
is not similarly small for the Grassmannian.

We now proceed to the proof of the main result in this section. For
this, we use a union bound to combine the probabilistic guarantees of
Lemma~\ref{lemma:lipshitz} with the estimates of the covering number
of $\Sigma(U,V)$.

\begin{proof}[of Theorem~\ref{thm:ripscaling}]

Let $\Omega =\{(V_i,W_i)\}$ be a finite set of subspaces that
satisfies the conditions of Lemma~\ref{lemma:covering} for
$\epsilon>0$, so $|\Omega| \leq \mathfrak{N}(\epsilon)$. For each
pair $(V_i,W_i)$, define the set of matrices
\begin{equation}
    \mathcal{B}_i:=
    \left\{X \,\, \Big|\,\, \exists (V,W) \,\, \mbox{such that}\,\,  X\in\Sigma(V,W)  \,\, \mbox{and}\,\, \rho(\Sigma(V,W),\Sigma(V_i,W_i))\leq \epsilon \right\}.
\end{equation}
Since $\Omega$ is an $\epsilon$-net, we have that the union of all
the $\mathcal{B}_i$ is equal to $\Sigma_{mnr}$.  Therefore, if for
all $i$, $(1-\delta)\|X\|_F \leq \|\cA(X)\|\leq (1+\delta)\|X\|_F$
for all $X\in\mathcal{B}_i$, we must have that $\delta_r(\cA)\leq
\delta$ proving that
\begin{equation}\label{eq:delta-prob}
\begin{split}
    \mathbf{P}(\delta_r(\cA) \leq \delta) &= \mathbf{P}\left[(1-\delta)\|X\|_F \leq \|\cA(X)\|\leq
(1+\delta)\|X\|_F \quad \forall \,\, X \mbox{ s.t. } \rank{X}\leq
    r \right]\\
    &\geq
    \mathbf{P}\left[\forall i\,\, (1-\delta)\|X\|_F \leq \|\cA(X)\|\leq
(1+\delta)\|X\|_F \quad\forall \,\, X \in \mathcal{B}_i\right]
    \end{split}
\end{equation}
Now note that if we have  $(1+ \|\cA\|)\epsilon \leq \delta/2$ and,
for all $Y \in \Sigma(V_i,W_i)$,  $(1-\delta/2)\|Y\|_F \leq
\|\cA(Y)\|\leq (1+\delta/2)\|Y\|_F$,
Lemma~\ref{lemma:lipshitz} implies that $(1-\delta)\|X\|_F \leq
\|\cA(X)\| \leq (1+\delta) \|X\|_F$ for all
$X\in\mathcal{B}_i$.  Therefore, using a union bound,
\eq{delta-prob} is greater than or equal to
\begin{equation}
\begin{split}
 1 &-
\sum_{i=1}^{|\Omega|} \mathbf{P}\left[\exists\,
    Y\in \Sigma(V_i,W_i)\quad
    \|\cA(Y)\| < (1-\frac{\delta}{2})\|Y\|_F\quad\mbox{or}\quad
    \|\cA(Y)\| > (1+\frac{\delta}{2})\|Y\|_F
     \right]\\ &\qquad\qquad-
    \mathbf{P}\left[\|\cA\| \geq
    \frac{\delta}{2\epsilon}-1\right]\, .
    \end{split}
\end{equation}
We can bound these quantities separately.  First we have by
Lemmas~\ref{lemma:lipshitz} and~\ref{lemma:covering}
\begin{equation}\label{eq:prob-eps-net}
\begin{split}
   \sum_{i=1}^{|\Omega|} \mathbf{P}&\left[\exists\,
    Y\in \Sigma(V_i,W_i)\quad
    \|\cA(Y)\| < (1-\frac{\delta}{2})\|Y\|_F\quad\mbox{or}\quad
    \|\cA(Y)\| > (1+\frac{\delta}{2})\|Y\|_F
     \right]\\
     &\leq 2 \mathfrak{N}(\epsilon)
    \left(\frac{24}{\delta}\right)^{r^2}
    \exp\left(-\frac{p}{2}(\delta^2/32 - \delta^3/96)\right)\\
    &\leq 2 \left(\frac{2C_0}{\epsilon}\right)^{r(m+n-2r)}
    \left(\frac{24}{\delta}\right)^{r^2}
    \exp\left(-\frac{p}{2}(\delta^2/32 -\delta^3/96)\right)\,.
    \end{split}
\end{equation}
Secondly, since $\cA$ is nearly isometric, there exists a
constant $\gamma$ such that
\begin{equation}
    \mathbf{P}\left( \|\cA\| \geq 1 + \sqrt{\frac{mn}{p}} + t \right) \leq \exp(-\gamma p
    t^2)\,.
\end{equation}
In particular,
\begin{equation}
    \mathbf{P}\left( \|\cA\| \geq \frac{\delta}{2\epsilon} -1 \right) \leq
    \exp\left(-\gamma p \left(\frac{\delta}{2\epsilon} - \sqrt{\frac{mn}{p}} - 2
    \right)^2\right)\,.
\end{equation}

We now must pick a suitable resolution $\epsilon$ to guarantee that
this probability is less than $\exp(-c_1 p)$ for a suitably chosen
constant $c_1$.  First note that if we choose $\epsilon<
(\delta/4)(\sqrt{mn/p}+1)^{-1}$,
\begin{equation}
     \mathbf{P}\left( \|\cA\| \geq \frac{\delta}{2\epsilon} -1 \right)
    \leq \exp(-\gamma mn)\, ,
\end{equation}
which achieves the desired scaling because $mn>p$. With this choice
of $\epsilon$, the quantity in Equation~\eq{prob-eps-net} is less
than or equal to
\begin{equation}
\begin{split}
&2 \left(\frac{8C_0 (\sqrt{mn/p}+1)}{\delta}\right)^{r(m+n-2r)}
     (24/\delta)^{r^2}\exp\left(-\frac{p}{2}(\delta^2/32 -
    \delta^3/96)\right)\\
&=\exp\left( -p a(\delta) +
r(m+n-2r)\log\left(\sqrt{\frac{mn}{p}}+1\right)\right.\\
&\qquad\qquad\qquad\qquad\qquad \left.+ r(m+n-2r)\log\left(\frac{
8C_0}{\delta}\right) +r^2\log\left(\frac{24}{\delta}\right)\right)
\end{split}
\end{equation}
where $a(\delta) = \delta^2/64-\delta^3/192$.  Since $mn/p<mn$ for
all $p>1$, there exists a constant $c_0$ independent of $m$,$n$,$p$,
and $r$, such that the sum of the last three terms in the exponent
are bounded above by $(c_0/a(\delta)) r(m+n)\log(mn)$. It follows
that there exists a constant $c_1$ independent of $m$,$n$,$p$, and
$r$ such that $p \geq c_0 r(m+n)\log(mn)$ observations are
sufficient to yield an RIP of $\delta$ with probability greater than
$1-e^{-c_1p}$.
\end{proof}

Heuristically, the scaling $p=O\left(r(m+n)\log(mn)\right)$ is very
reasonable, since a rank $r$ matrix has $r(m+n-r)$ degrees of
freedom. This coarse tail bound only provides asymptotic estimates
for recovery, and is quite conservative in practice. As we
demonstrate in Section \ref{sec:experiments}, minimum rank solutions
can be determined from between $2r(m+n-r)$ to $4r(m+n-r)$
observations for many practical problems.

\section{Algorithms for nuclear norm minimization}
\label{sec:algorithms}

A variety of methods can be developed for the effective minimization
of the nuclear norm over an affine subspace of matrices, and we do
not have room for a comprehensive treatment here. Instead, we focus
on three methods highlighting the trade-offs between computational
speed and guarantees on accuracy of the resulting solution. Directly
solving the semidefinite characterization of the nuclear norm
problem using primal-dual interior point methods is a numerically
efficient method for small problems and can be used to yield
accuracy up to floating-point precision.

Since interior point methods use second order information, the
memory requirements for computing descent directions quickly becomes
too large as the problem size increases. Moreover, for larger
problem sizes it is preferable to use methods that exploit, at least
partially, the structure of the problem. This can be done at several
levels, either by taking into account further information that may
be available about the linear map $\cA$ (e.g., the case of partially
observed Fourier measurements) or by formulating algorithms that are
specific to the nuclear norm problem. For the latter, we show how to
apply subgradient methods to minimize the nuclear norm over an
affine set. Such first-order methods cannot yield as high numerical
precision as interior point methods, but much larger problems can be
solved because no second-order information needs to be stored. For
even larger problems, we discuss a low-rank semidefinite programming
that explicitly works with a factorization of the decision variable.
This method can be applied even when the matrix decision variable
cannot fit into memory, but convergence guarantees are much less
satisfactory than in the other two cases.

\subsection{Interior Point Methods for Semidefinite programming}

For small problems where a high-degree of numerical precision is
required, interior point methods for semidefinite programming can be
directly applied to solve affine nuclear minimization problems. As
we have seen in earlier sections, the nuclear norm minimization
problem can be directly posed as a semidefinite programming problem
via the standard form primal-dual pair~\eq{nucnorm-sdp-pd}.  As
written, the primal problem has one $(n+m) \times (n+m)$
semidefinite constraint and $p$ affine constraints. Conversely, the
dual problem has one $(n+m) \times (n+m)$ semidefinite constraint
and $p$ scalar decision variables.  Thus, the total number of
decision variables (primal and dual) is equal to
$\binom{n+m+1}{2}+p$.

Modern interior point solvers for semidefinite programming generally
use primal-dual methods, and compute an update direction for the
current solution by solving a suitable Newton system. Depending on
the structure of the linear mapping $\cA$, this may entail
solving a potentially large, dense linear system.

If the matrix dimensions $n$ and $m$ are not too large, then any
good interior point SDP solver, such as SeDuMi \cite{sedumi} or
SDPT3 \cite{SDPT3}, will quickly produce accurate solutions.  In
fact, as we will see in the next section, problems with $n$ and $m$
around $50$ can be solved to machine precision in minutes on a
desktop computer to machine precision.  However, solving such a
primal-dual pair of programs with traditional interior point methods
can prove to be quite challenging when the dimensions of the matrix
$X$ are much bigger than $100 \times 100$, since in this case the
corresponding Newton systems become quite large. In the absence of
any specific additional structure, the memory requirements of such
dense systems quickly limit the size of problems that can be solved.

Perhaps the most important drawback of the direct SDP approach is
that it completely ignores the possibility of efficiently computing
the nuclear norm via a singular value decomposition, instead of the
less efficient eigenvalue decomposition of a bigger matrix. The
method we discuss next will circumvent this obstacle, by directly
working with subgradients of the nuclear norm.

\subsection{Projected subgradient methods}

The nuclear norm minimization \eq{nucnorm-prob} is a linearly
constrained nondifferentiable convex problem. There are numerous
techniques to approach this kind of problems, depending on the
specific nature of the constraints (e.g., dense vs. sparse), and the
possibility of using first- or second-order information.

In this section we describe a simple, easy to implement, subgradient
projection approach to the solution of \eq{nucnorm-prob}. This
first-order method will proceed by computing a sequence of feasible
points $\{X_k\}$, with iterates satisfying the update rule
\[
X_{k+1} = \Pi (X_k - s_k Y_k), \qquad Y_k \in \partial ||X_k||_*,
\]
where $\Pi$ is the orthogonal projection onto the affine subspace
defined by the linear constraints $\cA(X) = b$, and $s_k>0$
is a stepsize parameter. In other words, the method updates the
current iterate $X_k$ by taking a step along the direction of a
subgradient at the current point and then projecting back onto the
feasible set. Alternatively, since $X_k$ is feasible, we can rewrite
this as
\[
X_{k+1} = X_k - s_k \Pi_\cA Y_k,
\]
where $\Pi_\cA$ is the orthogonal projection onto the kernel
of $\cA$. Since the feasible set is an affine subspace,
there are several options for the projection $\Pi_\cA$. For
small problems, one can precompute it using, for example, a QR
decomposition of the matrix representation of $\cA$ and
store it. Alternatively, one can solve a least squares problem at each
step by iterative methods such as conjugate gradients.

The subgradient-based method described above is extremely simple to
implement, since only a subgradient evaluation is required at every
step. The computation of the subgradient can be done using the
formula given in \eq{subdiffnuc} earlier, thus requiring only a
singular value decomposition of the current point $X_k$.

A possible alternative here to the use of the SVD for the
subgradient computation is to directly focus on the ``angular''
factor of the polar decomposition of $X_k$, using for instance the
Newton-like methods developed by Gander in \cite{GanderPolar}.
Specifically, for a given matrix $X_k$, the Halley-like iteration
\[
X \rightarrow X (X'X+3I) (3 X'X+I)^{-1}
\]
converges globally and quadratically to the polar factor of $X$, and
thus yields an element of the subdifferential of the nuclear norm.
This iteration method (suitable scaled) can be faster than a direct
SVD computation, particularly if the singular values of the initial
matrix are close to 1. This could be appealing since presumably only
a very small number of iterations would be needed to update the
polar factor of $X_k$, although the nonsmoothness of the
subdifferential is bound to cause some additional difficulties.

Regarding convergence, for general nonsmooth problems, subgradient
methods do not guarantee a decrease of the cost function at every
iteration, even for arbitrarily small step sizes (see, e.g.,
\cite[\S 6.3.1]{BertsekasNLP}), unless the minimum-norm subgradient
is used. Instead, convergence is usually shown through the decrease
(for small stepsize) of the distance from the iterates $X_k$ to any
optimal point.  There are several possibilities for the choice of
stepsize $s_k$. The simplest choice that can guarantee convergence
is to use a diminishing stepsize with an infinite travel condition
(i.e., such that $\lim_{k \rightarrow \infty} s_k =0$ and
$\sum_{k>0} {s_k}$ diverging).

Often times, even the computation of a singular value decomposition
or Halley-like iteration can be too computationally expensive.  The
next section proposes a reduction of the size of the search space to
alleviate such demands.  We must give up guarantees of convergence
for this convenience, but this may be an acceptable trade-off for
very large-scale problems.

\subsection{Low-rank parametrization}

We now turn to a method that works with an explicit low-rank
factorization of $X$.  This algorithm not only requires less storage
capacity and computational overhead than the previous methods, but
for many problems does not even require one to be able to store the
decision variable $X$ in memory.  This is the case, for example, in
the matrix completion problem where $\cA(X)$ is a subset of the
entries of $X$.

Given observations of the form $\cA(X)=b$ of an $m\times n$ matrix
$X$ of rank $r$, a possible search algorithm to find a suitable $X$
would be to find a factorization $X = LR'$, where $L$ is an $m
\times r$ matrix and $R$ an $n \times r$ matrix, such that the
equality constraints are satisfied.  Since there are many possible
such factorizations, we could search for one where the matrices $L$
and $R$ have Frobenius norm as small as possible, i.e., the solution
of the optimization problem
\begin{equation}
\begin{split}
\min_{L,R} &\quad \tfrac{1}{2}(\|L\|_F^2+\|R\|_F^2)\\
\mbox{s.t.} &\quad  \cA(LR') = b.
\end{split}
\label{eq:nucnorm-nc}
\end{equation}
Even though the cost function is convex, the constraint is not. Such
a problem is a nonconvex quadratic program, and it is not evidently
easy to optimize. We show below that the minimization of the nuclear
norm subject to equality constraints is in fact equivalent to this
rather natural heuristic optimization, as long as $r$ is chosen to
be sufficiently larger than the rank of the optimum of the nuclear
norm problem.

\begin{lemma}
Assume $r \geq \rank{X_0}$. The nonconvex quadratic optimization
problem \eq{nucnorm-nc} is equivalent to the minimum nuclear norm
relaxation \eq{nucnorm-prob}.
\end{lemma}
\begin{proof}
Consider any feasible solution $(L,R)$ of \eq{nucnorm-nc}. Then,
defining $W_1 := LL'$, $W_2 := RR'$, and $X := LR'$ yields a
feasible solution of the primal SDP problem \eq{nucnorm-sdp-pd} that
achieves the same cost. Since the SDP formulation is equivalent to
the nuclear norm problem, we have that the optimal value of
\eq{nucnorm-nc} is always greater than or equal to the nuclear norm
heuristic.

For the converse, we can use an argument similar to the proof of
Proposition~\ref{prop:dualnorm}. From the SVD decomposition $X^* = U
\Sigma V'$ of the optimal solution of the nuclear norm relaxation
\eq{nucnorm-prob}, we can explicitly construct matrices $L:=U
\Sigma^\frac{1}{2}$ and $R:=V \Sigma^\frac{1}{2}$ for
\eq{nucnorm-nc} that yield exactly the same value of the objective.
\end{proof}

The main advantage of this reformulation is to substantially
decrease the number of primal decision variables from $nm$ to
$(n+m)r$.  For large problems, this is quite a significant reduction
that allows us to search for matrices of rank smaller than the order
of $100$, and $n+m$ in the hundreds of thousands on a desktop
computer. However, this problem is nonconvex and potentially subject
to local minima.  This is not as much of a problem as it could be,
for two reasons. First recall from Theorem~\ref{thm:reclowrank} that
if $\delta_{2r}(\cA)<1$, there is a unique $X^*$ with rank
at most $r$ such that $\cA(X^*)=b$. Since any local minima
$(L^*,R^*)$ of \eq{nucnorm-nc} is feasible, we would have
$X^*=L^*{R^*}'$ and we would have found the minimum rank solution.
Second, we now present an algorithm that is guaranteed to converge
to a local minima for a judiciously selected $r$.  We will also
provide a sufficient condition for when we can construct an optimal
solution of \eq{nucnorm-sdp-pd} from the solution computed by the
method of multipliers.

\paragraph{SDPLR and the method of multipliers}
For general semidefinite programming problems, Burer and Monteiro
have developed in~\cite{Burer03,Burer05} a nonlinear programming
approach that relies on a low-rank factorization of the matrix
decision variable.  We will adapt this idea to our problem, to
provide a first-order Lagrangian minimization algorithm that
efficiently finds a local minima of~\eq{nucnorm-nc}. As a
consequence of the work in~\cite{Burer05}, it will follow that for
values of $r$ larger than the rank of the true optimal solution, the
local minima of~\eq{nucnorm-nc} can be transformed into global
minima of~\eq{nucnorm-sdp-pd} under the identification $W_1=LL'$,
$W_2=RR'$ and $Y=LR'$.  We summarize below the details of this
approach.

The algorithm employed is called the \emph{method of multipliers}, a
standard approach for solving equality constrained optimization
problems~\cite{BertsekasLagrangeBook}.  The method of multipliers
works with an augmented Lagrangian for \eq{nucnorm-nc}
\begin{equation}\label{eq:lagrange}
\mathcal{L}_a(L,R; y,\sigma) := \textstyle{\frac{1}{2}}
(||L||^2_F+||R||^2_F) -y' (\cA(LR') -b) + \frac{\sigma}{2}
||\cA(LR') -b ||^2,
\end{equation}
where the $y_i$ are arbitrarily signed Lagrange multipliers and
$\sigma$ is a positive constant.  A somewhat similar algorithm was
proposed by Rennie et al in~\cite{Rennie05} in the collaborative
filtering.  In this work, the authors minimize $\mathcal{L}_a$ with
$\sigma$ fixed and $y=0$ to serve as a regularized algorithm for
matrix completion. Remarkably, by deterministically varying $\sigma$
and $y$, this method can be adapted into an algorithm for solving
linearly constrained nuclear-norm minimization.

In the method of multipliers, one alternately minimizes the
augmented Lagrangian with respect to the decision variables $L$ and
$R$, and then increases the value of the penalty coefficient
$\sigma$ and updates $y$.  The augmented Lagrangian can be minimized
using any local search technique, and the partial derivatives are
particularly simple to compute.  Let $\hat{y} :=
y-\sigma(\cA(LR')-b)$.  Then we have
\begin{equation*}
\begin{split}
   \nabla_L \mathcal{L}_a &= L - \cA^*(\hat y) R\\
   \nabla_R \mathcal{L}_a &= R - \cA^*(\hat y)' L.
   \end{split}
\end{equation*}
To calculate the gradients, we first compute the constraint
violations $\cA(LR')-b$, then form $\hat{y}$, and finally
use the above equations to compute the gradients.

As the number of iterations tends to infinity, only feasible points
will have finite values of $\mathcal{L}_a$, and for any feasible
point, $\mathcal{L}_a(L,R)$ is equal to the original cost function
$(\|L\|_F^2 + \|R\|_F^2)/2$. The method terminates when $L$ and $R$
are feasible, as in this case the Lagrangian is stationary and we
are at a local minima of \eq{nucnorm-nc}.  Including the $y$
multipliers improves the conditioning of each subproblem where
$\mathcal{L}_a$ is minimized and enhances the rate of convergence.
The following theorem shows that when the method of multipliers
converges, it converges to a local minimum of~\eq{nucnorm-nc}.

\begin{theorem}\label{thm:converge}
Suppose we have a sequence $(L^{(k)},R^{(k)},y^{(k)})$ of local
minima of the augmented Lagrangian at each step of the method of
multipliers.  Assume that our sequence of $\sigma^{(k)}\rightarrow
\infty$ and that the sequence of $y^{(k)}$ is bounded. If
$(L^{(k)},R^{(k)})$ converges to $(L^*,R^*)$ and the linear map
\begin{equation}
 \Lambda^{(k)}(y):= \left[\begin{array}{c} \cA^*(y) R^{(k)}\\
\cA^*(y)' L^{(k)}\end{array}\right]
\end{equation}
has kernel equal to the zero vector for all $k$, then there exists a
vector $y^*$ such that
\begin{enumerate}\renewcommand{\labelenumi}{(\roman{enumi})}
\item\label{item:stationary} $\nabla \mathcal{L}_a(L^{*},R^{*};y^{*}) = 0$
\item $\cA(L^*{R^*}') = b$
\end{enumerate}
\end{theorem}
\begin{proof}
This proof is standard and follows the approach
in~\cite{BertsekasLagrangeBook}.  As above, we define
$\hat{y}^{(k)}:=y^{(k)} -
\sigma^{(k)}(\cA(L^{(k)}{R^{(k)}}')-b)$ for all $k$. Since
$(L^{(k)},R^{(k)})$ minimize the augmented Lagrangian at iteration
$k$, we have
\begin{equation}\label{eq:stationary-iteration}
\begin{split}
   0 &= L^{(k)} - \cA^*(\hat{y}^{(k)})R^{(k)}\\
   0 &= R^{(k)} - \cA^*(\hat{y}^{(k)})'L^{(k)},
   \end{split}
\end{equation}
which we may rewrite as
\begin{equation}
    \Lambda^{(k)}(\hat{y}^{(k)}) =
    \left[\begin{array}{c} L^{(k)} \\ R^{(k)} \end{array}\right]\,.
\end{equation}
Since we have assumed that there is no non-zero $y$ with
$\Lambda^{(k)}(y)= 0$, there exists a left-inverse and we can solve
for $\hat{y}^{(k)}$.
\begin{equation}
    \hat{y}^{(k)} =
    {\Lambda^{(k)}}^\dagger \left(\left[\begin{array}{c}
L^{(k)}\\R^{(k)}\end{array}\right]\right)\,.
\end{equation}
Everything on the right-hand side is bounded, and $L^{(k)}$ and
$R^{(k)}$ converge.  Therefore, we must have that $\hat{y}^{(k)}$
converges to some $\hat{y}^{*}$.  Taking the limit of
\eq{stationary-iteration} proves (i). To prove (ii), note that we
must have ${\hat{y}^{(k)}}$ is bounded. Since $y^{(k)}$ is also
bounded, we find that
$\sigma^{(k)}(\cA(L^{(k)}{R^{(k)}}')-b)$ is also bounded.
But $\sigma^{(k)}\rightarrow \infty$ implies that
$\cA(L^*{R^*}')=b$, completing the proof.
\end{proof}

Suppose the decision variables are chosen to be of size $m\times
r_d$ and $n\times r_d$. A necessary condition for $\Lambda^{k}(y)$
to be full rank is for the number of decision variables $r_d(m+n)$
to be greater than the number of equalities $p$.  In particular,
this means that we must choose $r_d \geq p/(m+n)$ in order to have
any hopes of satisfying the conditions of
Theorem~\ref{thm:converge}.

We close this section by relating the solution found by the method
of multipliers to the optimal solution of the nuclear norm
minimization problem.  We have already shown that when the low-rank
algorithm converges, it converges to a low-rank solution of
$\cA(X)=b$. If we additionally find that
$\cA^*(y^*)$ has norm less than or equal to one, then it is
dual feasible.  One can check using straightforward algebra that
$(L^*{R^*}',L^*{L^*}',R^*{R^*}')$ and $y^*$ form an optimal
primal-dual pair for~\eq{nucnorm-sdp-pd}.  This analysis proves the
following theorem.

\begin{theorem}
Let $(L^*,R^*,y^*)$ satisfy (i)-(ii) in Theorem~\ref{thm:converge}
and suppose $\|\cA^*(y^*)\|\leq 1$. Then
$(L^*{R^*}',L^*{L^*}',R^*{R^*}')$ is an optimal primal solution and
$y^*$ is an optimal dual solution of \eq{nucnorm-sdp-pd}.
\end{theorem}

\section{Numerical Experiments}
\label{sec:experiments}

\begin{figure}
\centering
\includegraphics[width=8cm]{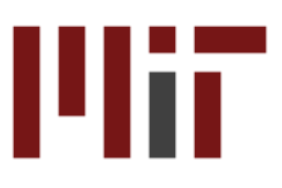}
\caption{\small The MIT logo image.  The associated matrix has
dimensions $46 \times 81$ and has rank $5$.}\label{fig:mit}
\end{figure}

To illustrate the scaling of low-rank recovery for a particular
matrix $M$, consider the MIT logo presented in Figure~\ref{fig:mit}.
The image has a total of $46$ rows and $81$ columns (total 3726
elements), and 3 distinct non-zero numerical values corresponding to
the colors white, red, and grey. Since the logo only has $5$
distinct rows, it has rank $5$. For each of the ensembles discussed
in Section~\ref{sec:gaussian}, we sampled measurement matrices with
$p$ ranging between $700$ and $1500$, and solved the semidefinite
program~\eq{sdp-embedding} using the freely available software
SeDuMi~\cite{sedumi}.  On a 2.0 GHz Laptop, each semidefinite
program could be solved in less than four minutes.  We chose to use
this interior point method because it yielded the highest accuracy
in the shortest amount of time, and we were interested in
characterizing precisely when the nuclear norm heuristic succeeded
and failed.

\begin{figure}
  \centering
  \begin{tabular}{c}
    \includegraphics[width=12cm]{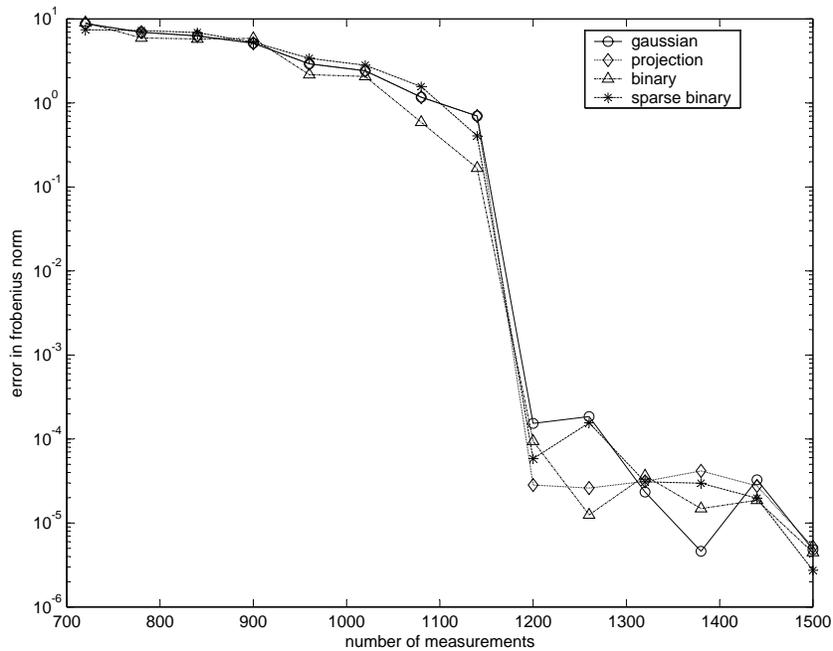}\\(a)\\
    \includegraphics[width=12cm]{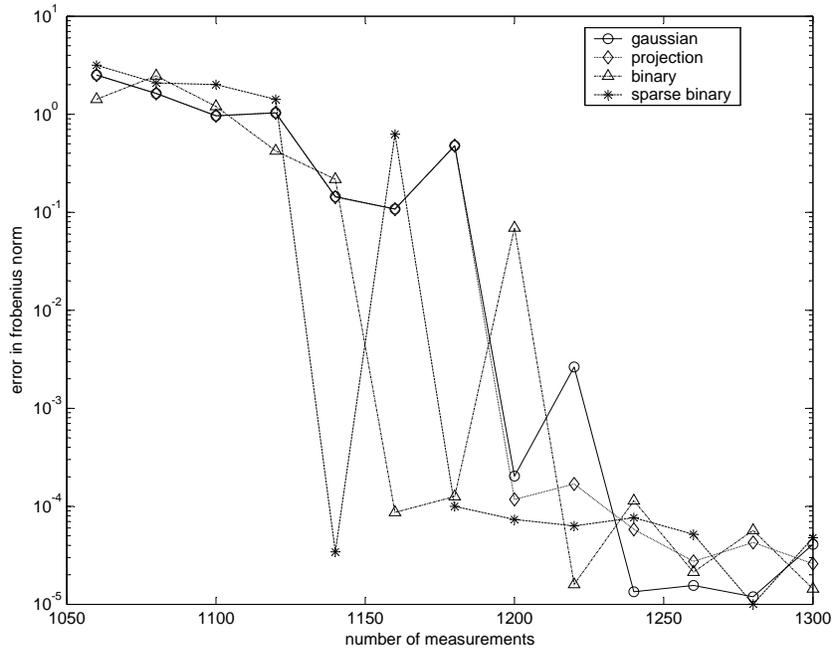}\\(b)\\
      \end{tabular}
  \caption{\small
(a) Error, as measured by the Frobenius norm, between the recovered
image and the ground truth.  Observe that there is a sharp
transition to near zero error at around $1200$ measurements.  (b)
Zooming in on this transition, we see fluctuation between high and
low error when between $1125$ and $1225$ measurements are
available.}
  \label{fig:mit-error}
\end{figure}

\begin{figure}
  \centering
  \begin{tabular}{ccc}
    \includegraphics[width=5cm]{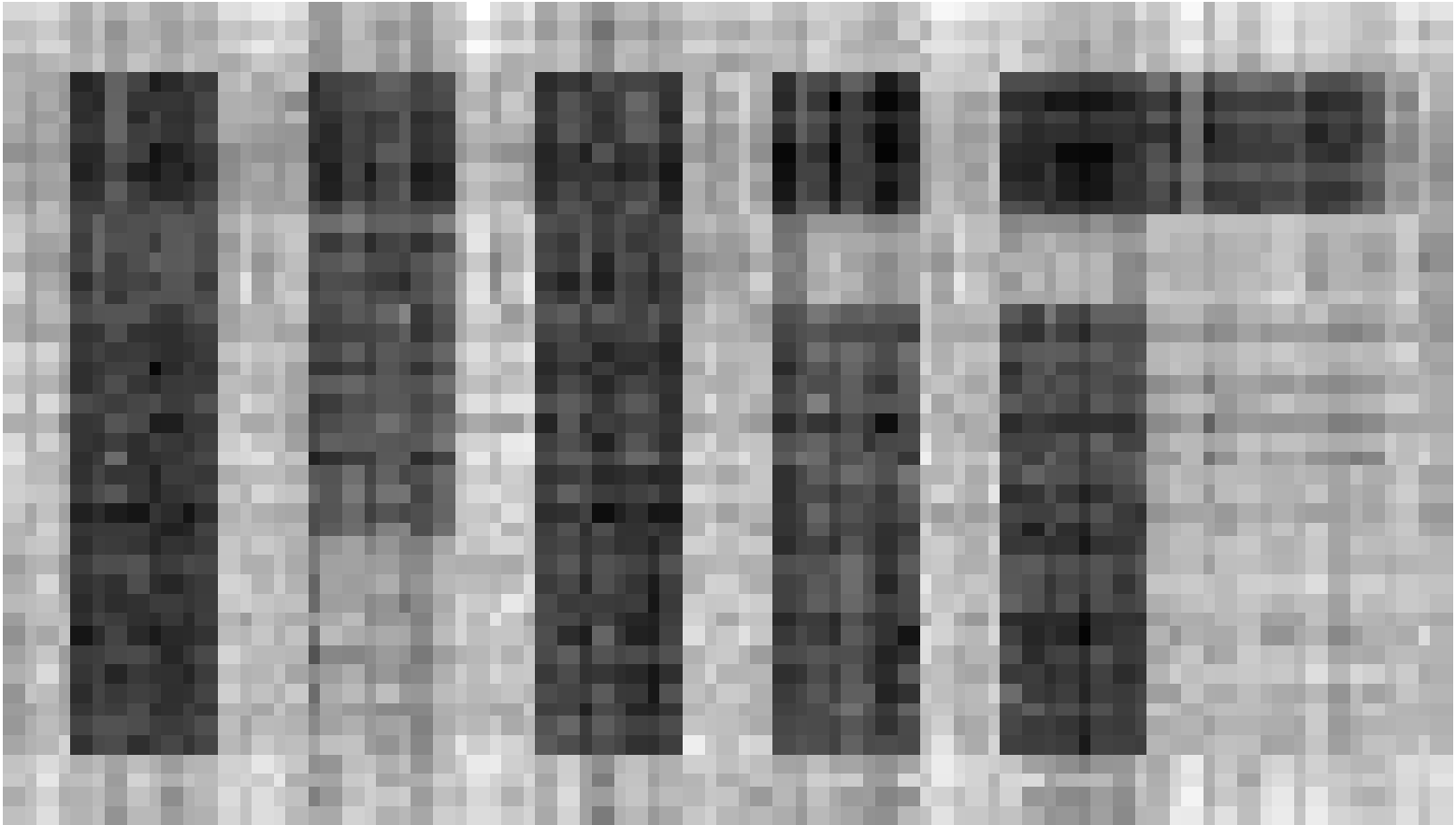} &
    \includegraphics[width=5cm]{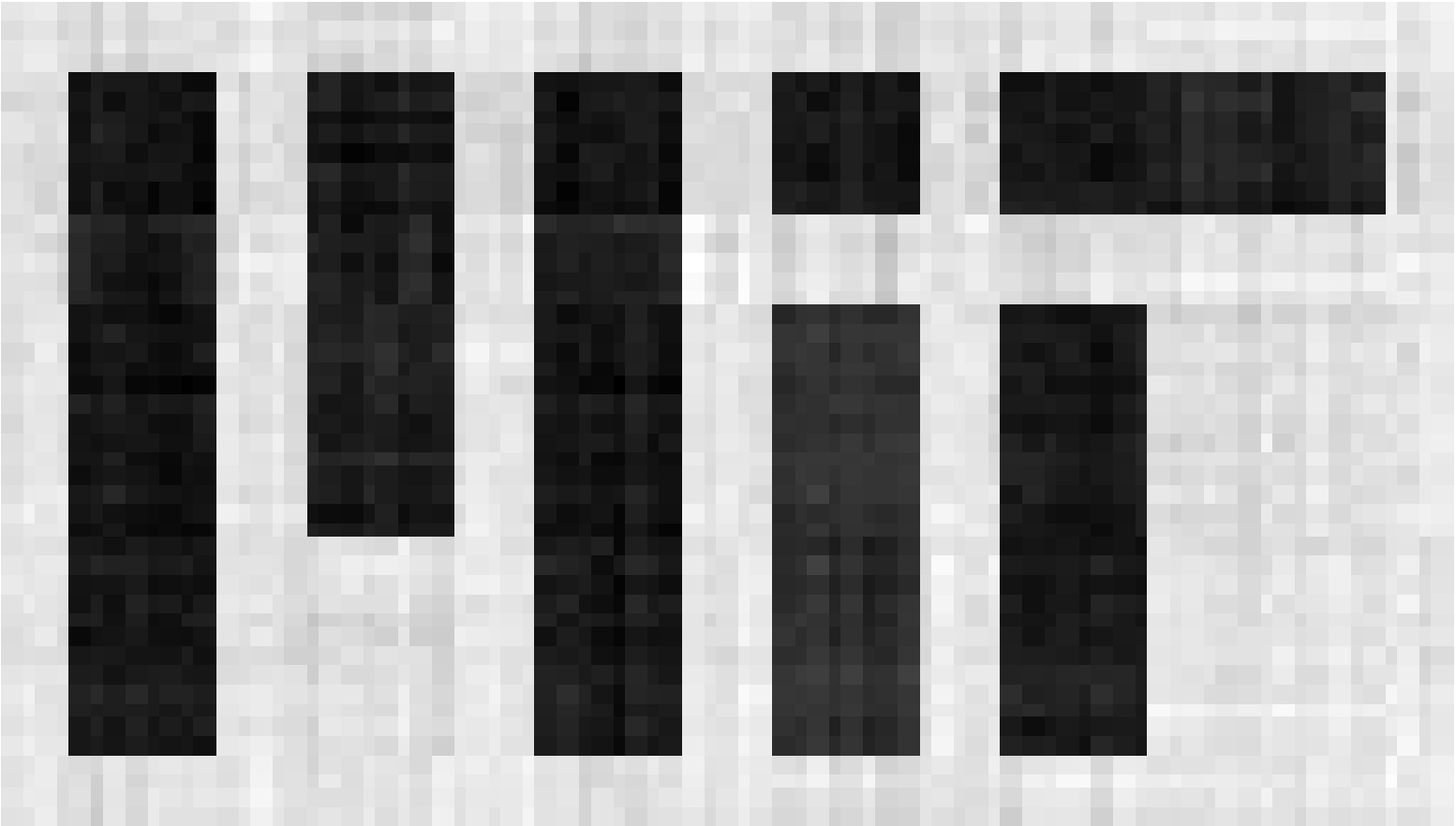} &
    \includegraphics[width=5cm]{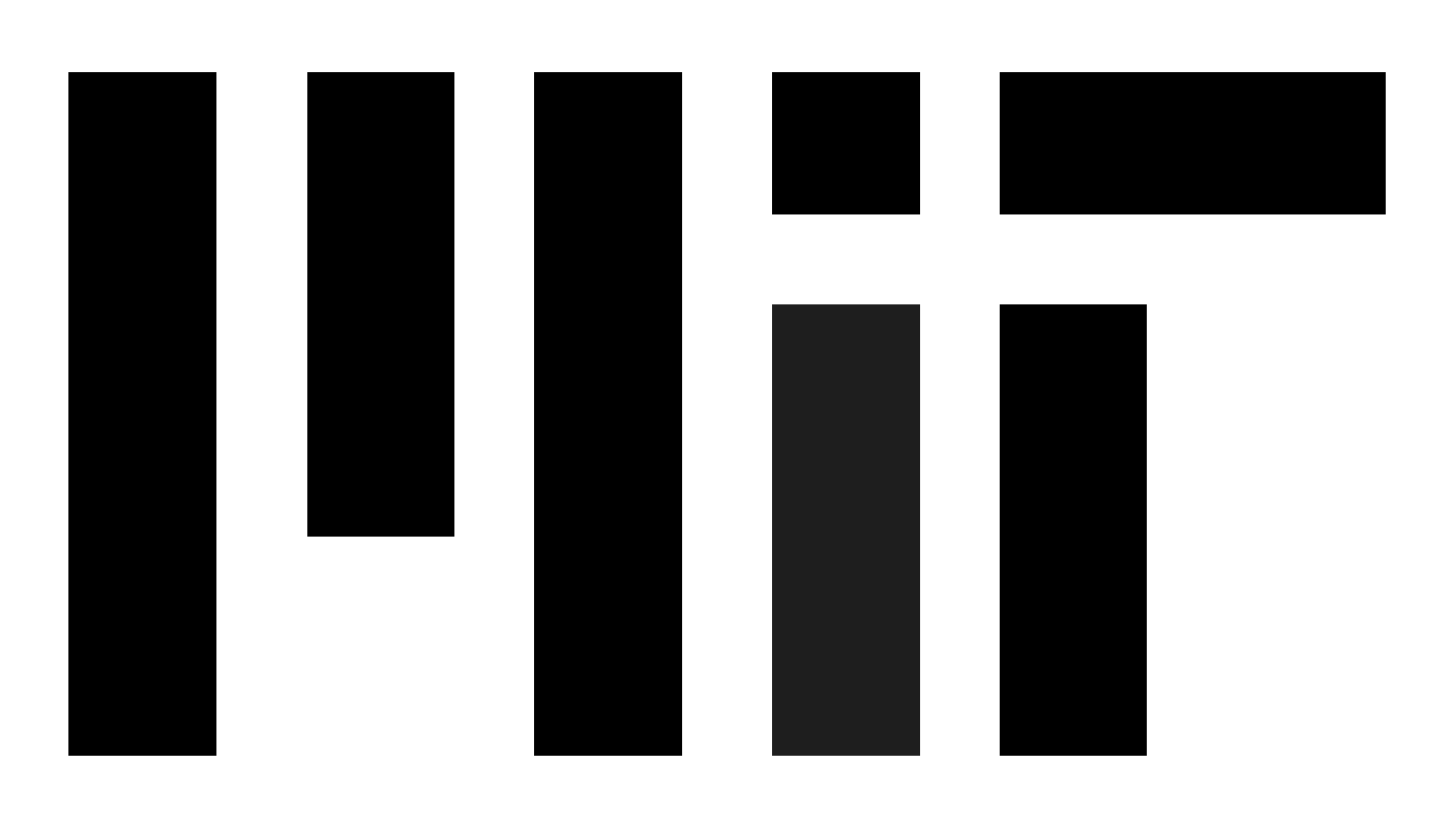}\\
    (a) & (b) & (c)\\
      \end{tabular}
  \caption{\small Example recovered images using the Gaussian
ensemble. (a) 700 measurements. (b) 1100 measurements (c) 1250
measurements. The total number of pixels is $46 \times 81 = 3726$.
Note that the error is plotted on a logarithmic scale.}
  \label{fig:mit-recovery}
\end{figure}

Figure~\ref{fig:mit-error} plots the Frobenius norm of the
difference between the optimal point of the semidefinite program and
the true image in Figure~\ref{fig:mit-error}.  We observe a sharp
transition to perfect recovery near 1200 measurements which is
approximately equal to $2r(m+n-r)$. In
Figure~\ref{fig:mit-recovery}, we graphically plot the recovered
solutions for various values of $p$ under the Gaussian ensemble.

To demonstrate the average behavior of low-rank recovery, we
conducted a series of experiments for a variety of the matrix sizes
$n$, ranks $r$, and numbers of measurements $p$.  For a fixed $n$,
we constructed random recovery scenarios for low-rank $n\times n$
matrices.  For each $n$, we varied $p$ between $0$ and $n^2$ where
the matrix is completely discovered.  For a fixed $n$ and $p$, we
generated all possible ranks such that $r(2n-r)\leq p$. This cutoff
was chosen because beyond that point there would be an infinite set
of matrices of rank $r$ satisfying the $p$ equations.

\begin{figure}
  \centering
  \begin{tabular}{c}
    \includegraphics[width=11cm]{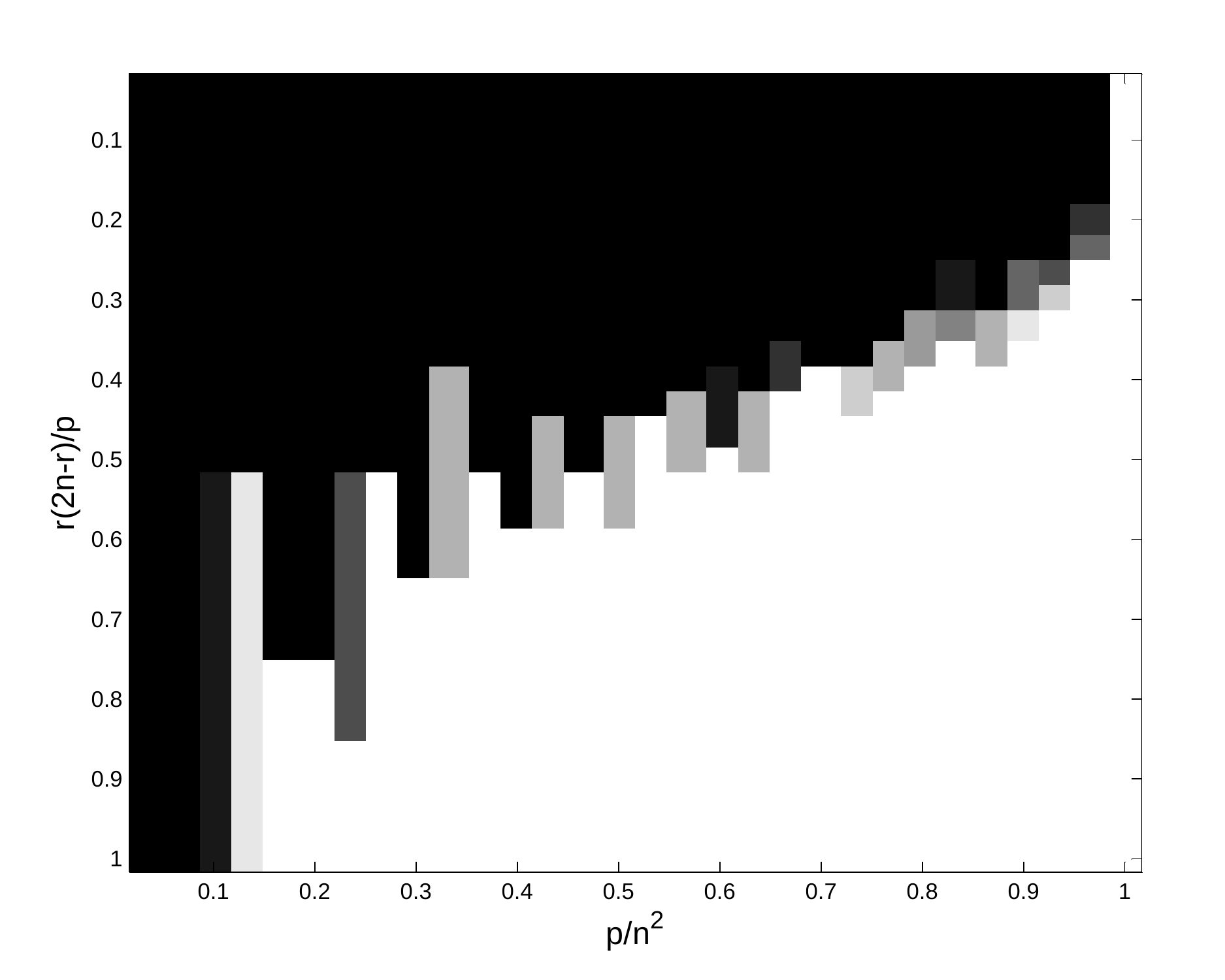}\\(a)\\
    \includegraphics[width=11cm]{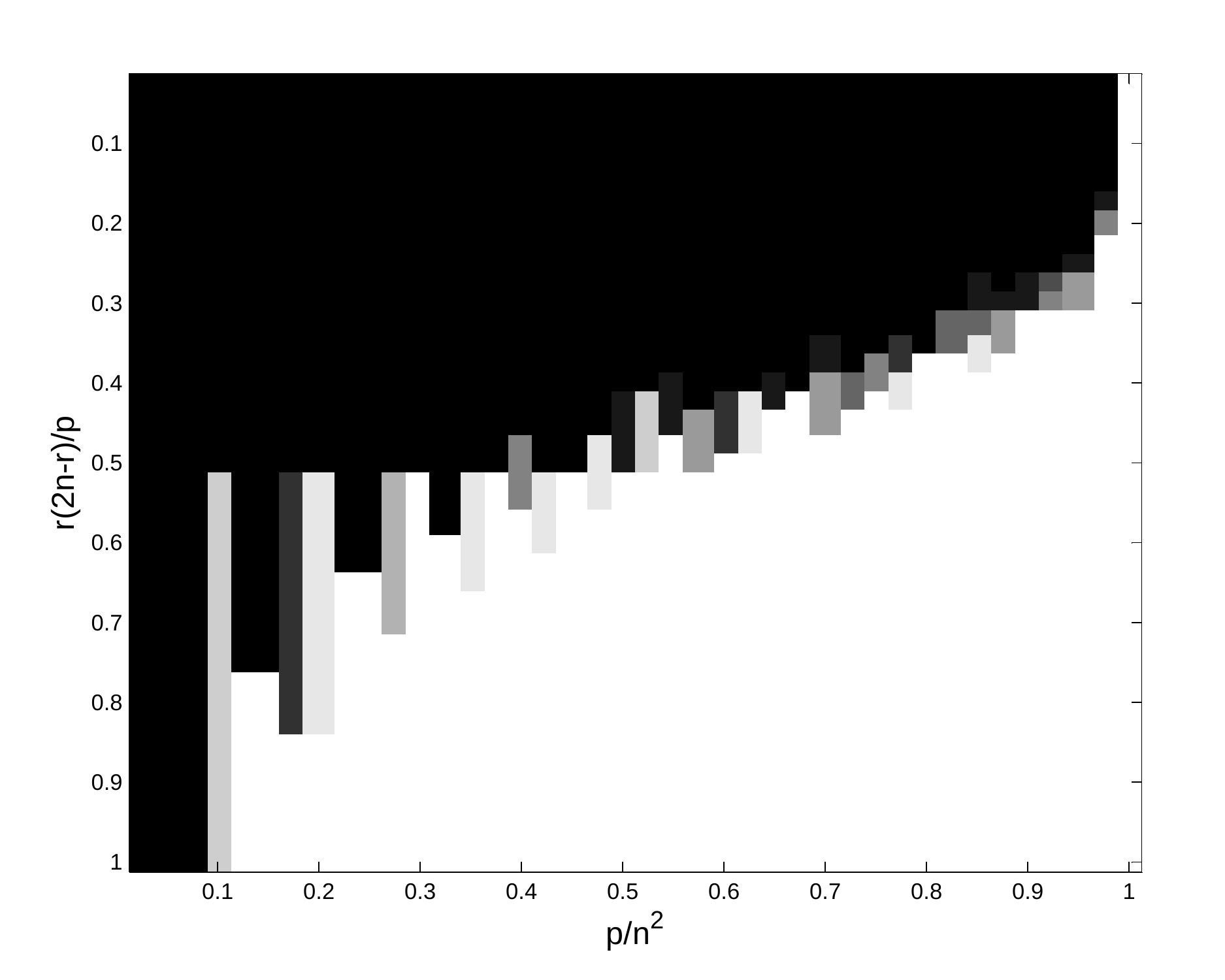}\\(b)\\
      \end{tabular}
  \caption{\small
For each $(n,p,r)$ triple, we repeated the following procedure ten
times. A matrix of rank $r$ was generated by choosing two random
$n\times r$ factors $Y_L$ and $Y_R$ with i.i.d. random entries and
set $Y_0=Y_LY_R'$. We select a matrix $\bfA$ from the Gaussian
ensemble with $p$ rows and $n^2$ columns. Then we solve the nuclear
norm minimization subject to $\bfA\vvec{X} = \bfA\vvec{Y_0}$ We
declare $Y_0$ to be recovered if $\|X-Y_0\|_F/\|Y_0\|_F<10^{-3}$.
The results are shown for (a) $n=30$ and (b) $n=40$. The color of
each cell reflects the empirical recovery rate (scaled between $0$
and $1$).  White denotes perfect recovery in all experiments, and
black denotes failure for all experiments.} \label{fig:phase-trans}
\end{figure}

For each $(n,p,r)$ triple, we repeated the following procedure $10$
times. A matrix of rank $r$ was generated by choosing two random
$n\times r$ factors $Y_L$ and $Y_R$ with i.i.d. random entries and
setting $Y_0=Y_LY_R'$. A matrix $\bfA$ was sampled from the Gaussian
ensemble with $p$ rows and $n^2$ columns.  Then the nuclear norm
minimization
\begin{equation}
    \begin{split}
        \min_{X} &\quad \|X\|_*\\
        \mbox{s.t.} &\quad \bfA\vvec{X} = \bfA\vvec{Y_0}
    \end{split}
\end{equation}
was solved using the SDP solver SeDuMi on the formulation
\eq{sdp-embedding}. Again, we chose to use SeDuMi because we wanted
to precisely distinguish between success and failure of the
heuristic.  We declared $Y_0$ to be recovered if
$\|X-Y_0\|_F/\|Y_0\|_F<10^{-3}$. Figure~\ref{fig:phase-trans} shows
the results of these experiments for $n=30$ and $40$.  The color of
the cell in the figures reflects the empirical recovery rate of the
$10$ runs (scaled between $0$ and $1$).  White denotes perfect
recovery in all experiments, and black denotes failure for all
experiments.

These experiments demonstrate that the logarithmic factors and
constants present in our scaling results are somewhat conservative.
For example, as one might expect, low-rank matrices are perfectly
recovered by nuclear norm minimization when $p=n^2$ as the matrix is
uniquely determined. Moreover, as $p$ is reduced slightly away from
this value, low-rank matrices are still recovered $100$ percent of
the time for most values of $r$.  Finally, we note that despite the
asymptotic nature of our analysis, our experiments demonstrate
excellent performance with low-rank matrices of size $30\times 30$
and $40 \times 40$ matrices, showing that the heuristic is practical
even in low-dimensional settings.

Intriguingly, Figure~\ref{fig:phase-trans} also demonstrates a
``phase transition'' between perfect recovery and failure.  As
observed in several recent papers by Donoho and his collaborators
(See e.g.~\cite{DT1,DT2}), the random sparsity recovery problem has
two distinct connected regions of parameter space: one where the
sparsity pattern is perfectly recovered, and one where no sparse
solution is found.  Not surprisingly, Figure~\ref{fig:phase-trans}
illustrates an analogous phenomenon in rank recovery.  Computing
explicit formulas for the transition between perfect recovery and
failure is left for future work.

\section{Discussion and future developments}
\label{sec:conclusions}

Having illustrated the natural connections between affine rank
minimization and affine cardinality minimization, we were able to
draw on these parallels to determine scenarios where the nuclear
norm heuristic was able to exactly solve the rank minimization
problem.  These scenarios directly generalized conditions for which
the $\ell_1$ heuristic succeeded and ensembles of linear maps for
which these conditions hold. Furthermore, our experimental results
display similar recovery properties to those demonstrated in the
empirical studies of $\ell_1$ minimization. Inspired by the success
of this program, we close this report by briefly discussing several
exciting directions that are natural continuations of this work
building on more analogies from the compressed sensing literature.
We also describe possible extensions to more general notions of
parsimony.

\paragraph{Factored measurements and alternative ensembles}
All of the measurement ensembles require the storage of $O(mnp)$
numbers.  For large problems this is wholly impractical.  There are
many promising alternative measurement ensembles that seem to obey
the same scaling laws as those presented in
Section~\ref{sec:gaussian}. For example, ``factored'' measurements,
of the form $A_i:X \mapsto u_i^T X v_i$, where $u_i, v_i$ are
Gaussian random vectors empirically yield the same performance as
the Gaussian ensemble. This factored ensemble only requires storage
of $O((m+n)p)$ numbers, which is a rather significant savings for
very large problems.  The proof in Section~\ref{sec:gaussian} does
not seem to extend to this ensemble, thus new machinery must be
developed to guarantee properties about such low-rank measurements.

\paragraph{Noisy Measurements and Low-rank approximation}
Our results in this paper address only the case of exact (noiseless)
measurements. It is of natural interest to understand the behavior
of the nuclear norm heuristic in the case of noisy data. Based on
the existing results for the sparse case (e.g., \cite{Candes05b}),
it would be natural to expect similar stability properties of the
recovered solution, for instance in terms of the $\ell_2$ norm of
the computed solution.  Such an analysis could also be used to study
the nuclear norm heuristic as an approximation technique where a
matrix has rapidly decaying singular values and a low-rank
approximation is desired.

\paragraph{Incoherent Ensembles and Partially Observed Transforms}
Again, taking our lead from the compressed sensing literature, it
would be of great interest to extend the results of~\cite{Romberg07}
to low-rank recovery.  In this work, the authors show that partially
observed unitary transformations of sparse vectors can be used to
recover the sparse vector using $\ell_1$ minimization.  There are
many practical applications where low-rank processes are partially
observed.  For instance, the matrix completion problem can be
thought of as partial observations under the identity
transformations.  As another example, there are many examples in
two-dimensional Fourier spectroscopy where only partial information
can be observed due to experimental constraints.

\paragraph{Alternative numerical methods}
Besides the techniques described in Section~\ref{sec:algorithms},
there are a number of interesting additional possibilities to solve
the nuclear norm minimization problem. An appealing suggestion is to
combine the strength of second-order methods (as in the SDP
approach) with the known geometry of the nuclear norm (as in the
subgradient approach), and develop a customized interior point
method, possibly yielding faster convergence rates, while still
being relatively memory-efficient.

It is also of much interest to investigate the possible adaptation
of some of the successful path-following approaches in traditional
$\ell_1$/cardinality minimization, such as the Homotopy
\cite{Osborne00} or LARS (least angle regression) \cite{Efron04}.
This may be not be completely straightforward, since the efficiency
of many of these methods often relies explicitly on the polyhedral
structure of the feasible set of the $\ell_1$ norm problem.

\paragraph{Geometric interpretations}
For the case of cardinality/$\ell_1$ minimization, a beautiful
geometric interpretation has been set forth by Donoho and Tanner
\cite{DT1,DT2}. Key to their results is the notion of \emph{central
$k$-neighborliness} of a centrosymmetric polytope, namely the
property that every subset of $k+1$ vertices not including an
antipodal pair spans a $k$-face. In particular, they show that the
$\ell_1$ heuristic always succeeds whenever the image of the
$\ell_1$ unit ball (the cross-polytope) under the linear mapping
$\cA$ is a centrally $k$-neighborly polytope.

In the case of rank minimization, the direct application of these
concepts fails, since the unit ball of the nuclear norm is not a
polyhedral set. Nevertheless, it seems likely that a similar
explanation could be developed, where the key feature would be the
preservation under a linear map of the extremality of the components
of the boundary of the nuclear norm unit ball defined by low-rank
conditions.

\paragraph{Jordan algebras}
As we have seen, our results for the rank minimization problem
closely parallel the earlier developments in cardinality
minimization. A convenient mathematical framework that allows the
simultaneous consideration of these cases as well as a few new ones,
is that of \emph{Jordan algebras} and the related symmetric cones
\cite{FarautKoranyi}. In the Jordan-algebraic setting, there is an
intrinsic notion of rank that agrees with the cardinality of the
support in the case of the nonnegative orthant or the rank of a
matrix in the case of the positive semidefinite cone. Besides
mathematical elegance, a direct Jordan-algebraic approach would
transparently yield similar results for the case of second-order (or
Lorentz) cone constraints.

As specific examples of the power and elegance of this approach, we
mention the work of Faybusovich \cite{Faybu97} and Schmieta and
Alizadeh \cite{SchmietaAlizadeh} that provide a unified development
of interior point methods for symmetric cones, as well as
Faybusovich's work on convexity theorems for quadratic mappings
\cite{FaybuQM}.

\paragraph{Parsimonious models and optimization}
Sparsity and low-rank are two specific classes of parsimonious (or
low-complexity) descriptions. Are there other kinds of
easy-to-describe parametric models that are amenable to exact
solutions via convex optimizations techniques? Given the intimate
connections between linear and semidefinite programming and the
Jordan algebraic approaches described earlier, it is likely that
this will require alternative tractable convex optimization
formulations.

\section{Acknowledgements}
We thank Stephen Boyd, Emmanuel Cand\`es, Jos\'{e} Costa, John
Doyle, Ali Jadbabaie, Ali Rahimi, and Michael Wakin for their useful
comments and suggestions. We also thank the IMA in Minneapolis for
hosting us during the initial stages of our collaboration.

{\small
\bibliographystyle{abbrv}
\bibliography{brecht}
}
\end{document}